\documentclass[preprint,12pt,3p]{elsarticle}

\def\x{{\bm x}}
\def\hx{{\hat{\bm x}}}
\def\bxi{{\bm \xi}}
\def\d{{\mbox{d}}}
\def\u{{\bm u}}

\def\hf{{\hat{f}}}
\def\hx{{\hat{\x}}}
\def\hc{{\hat{c}}}
\def\hht{{\hat{t}}}
\def\hQ{{\hat{Q}}}
\def\hT{{\hat{T}}}
\def\hu{{\hat{\u}}}
\def\hxi{{\hat{\bxi}}}
\def\htau{{\hat{\tau}}}
\def\hrho{{\hat{\rho}}}

\usepackage{amsmath,amssymb}
\usepackage{lipsum,soul,ulem}
\usepackage{bm}
\usepackage{graphicx}
\usepackage{graphics}
\usepackage{color}

\usepackage{lipsum}
\usepackage{amsfonts}
\usepackage{epstopdf}
\usepackage{algorithmic}

\newtheorem{theorem}{\textbf{Theorem}}
\newtheorem{definition}{\textbf{Definition}}
\newtheorem{remark}{\textbf{Remark}}
\newproof{proof}{\textbf{Proof}}
\journal{Journal of Computational Physics}

\begin{document}

\begin{frontmatter}
	
\title{On unified preserving properties of kinetic schemes}
\author[add1]{Zhaoli Guo}
\ead{zlguo@hust.edu.cn}
\author[add2]{Jiequan Li}
\ead{li\_jiequan@iapcm.ac.cn}
\author[add3]{Kun Xu}
\ead{makxu@ust.hk}
\address[add1]{State Key Laboratory of Coal Combustion, School of Energy and Power Engineering, Huazhong University of Science and Technology,Wuhan, 430074, China}
\address[add2]{Laboratory of Computational Physics, Institute of Applied Physics and Computational Mathematics, Beijing, 100088 China; and Center for Applied Physics and Technology, Peking University, 100871, China}
\address[add3]{Department of Mathematics, Hong Kong University of Science and Technology, Clear Water Bay, Hong Kong, China}

\date{\today}

\begin{abstract}
 The kinetic theory provides a physical basis for developing multiscal methods for gas flows covering a wide range of flow regimes. A particular challenge for kinetic schemes is whether they can capture the correct hydrodynamic behaviors of the system in the continuum regime (i.e., as the Knudsen number $\epsilon\ll 1$ ) without enforcing kinetic scale resolution. At the current stage, {the main approach to analyze such a property is the asymptotic preserving (AP) concept, which aims to show whether a kinetic scheme reduces to a solver for the hydrodynamic equations as $\epsilon \to 0$, such as the shock capturing scheme for the Euler equations. However, the detailed asymptotic properties of the kinetic scheme are indistinguishable when $\epsilon$ is small but finite under the AP framework}. In order to distinguish different characteristics of kinetic schemes, in this paper we introduce the concept of unified preserving (UP) aiming at assessing asmyptotic orders of a kinetic scheme by employing the modified equation approach and Chapman-Enskon analysis. It is shown that the UP properties of a kinetic scheme generally depend on the spatial/temporal accuracy and closely on the inter-connections among the three scales (kinetic scale, numerical scale, and hydrodynamic scale) and their corresponding coupled dynamics. Specifically, the numerical resolution and specific discretization of particle transport and collision determine the flow physics of the scheme in different regimes, especially in the near continuum limit. As two examples, the UP methodology is applied to analyze the discrete unified gas-kinetic scheme and a second-order implicit-explicit Runge-Kutta scheme in their asymptotic behaviors in the continuum limit.
\end{abstract}


\end{frontmatter}

\section{Introduction}
In recent years there are increasing interests on simulating multiscale gas flows in different flow regimes covering a wide range of Knudsen numbers,{defined as $\epsilon=\hat{\lambda}_0 / \hat{l}_0 \sim \hat{\tau}_0 /\hat{t}_0$, where $\hat{\lambda}_0$ and $\hat{\tau}_0$ are the typical mean free path and collision time of gas molecules, $\hat{l}_0$ and $\hat{t}_0$ are the characteristic hydrodynamic length and time scales, respectively.} It is a challenging problem for modeling and simulating such flows due to the large spans of temporal and/or spatial scales variations under a wide range of physical scenarios, and it becomes even more complicated at discrete levels with involvement of mesh size and time step scales \cite{ref:Xu-POF2017}. The classical computational fluid dynamics methods based on the Euler or Navier-Stokes equations are limited to continuum flows, while the kinetic particle methods, such as the Direct Simulation Monte Carlo, are mainly suitable for rarefied flows and encounter difficulties for continuum and near continuum flows.

On the other hand, it is well understood that gas kinetic models (Boltzmann or model equations) defined on the kinetic scale can lead to the Euler and Navier-Stokes equations on the hydrodynamic scale, and the gas kinetic theory provides a solid basis for developing schemes uniformly for flows from kinetic to hydrodynamic regimes. Therefore, developing deterministic numerical methods based on gas kinetic theory, which we call {\it kinetic schemes}, have attracted much attention. Actually, a variety of such schemes have been developed from different point of views in recent years, such as the lattice Boltzmann equation (LBE) method \cite{ref:GuoBook}, gas-kinetic scheme (GKS) \cite{ref:GKS}, semi-Lagrangian method \cite{ref:semi-Lagrangian}, Implicit-Explicit (IMEX) method \cite{ref:IMEX13,ref:IMEX17,ref:IMEX07,ref:IMEX_Hu18}, unified gas-kinetic scheme (UGKS)\cite{ref:UGKS}, and discrete unified gas-kinetic scheme (DUGKS) \cite{ref:DUGKS15,ref:DUGKS13,ref:DUGKS_Rev}. It is noted that the LBE and GKS are kinetic solvers for the Navier-Stokes equations, while the UGKS and DUGKS are designed mainly for multiscale flows. The progress of numerical methods based on the kinetic equations can be found in a recent review article \cite{ref:Review14} and references therein.

{
The asymptotic behavior of the Boltzmann equation at small $\epsilon$ can be assessed by asymptotic analysis under certain scaling of the flow, such as the Chapman-Enskog analysis under acoustic scaling. Generally, the Chapman-Enskog solution of the non-dimensional Boltzmann equation can be expressed as $f(\epsilon; \x, t)=\sum_{k=0}^{\infty}{\epsilon^k f^{(k)}(\x, t)}$, where the coefficients $f^{(k)}$ can be determined iteratively at the consecutive order of $\epsilon$. Here $\x=\hat{\x}/\hat{l}_0$ and $t=\hat{t}/\hat{t}_0$ are the dimensionless space and time variables, respectively, while the hatted variables are the corresponding dimensional ones. Although the Chapman-Enskog approximation may deviate from the true solution, it can provide some useful asymptotics information at small $\epsilon$. For instance, the solution is consistent with the Euler and Navier-Stokes equations at the zeroth and first orders of $\epsilon$, respectively.}

For a kinetic scheme, {in addition to the kinetic and hydrodynamic scales, numerical time and length scales (time step $\Delta \hht$ and mesh size $\Delta \hx$) are also involved, and its asymptotic behaviors will become more complicated. Furthermore, its capability of capturing accurate solutions in the transition regime between kinetic and hydrodynamic ones depends closely on its limiting solution in the continuum regime with the hydrodynamic scale resolution. Therefore, it is critical to make a reliable assessment of the asymptotic behaviors of a kinetic scheme.}
{Generally, a numerical solution $f_h$ of the (non-dimensional) Boltzmann equation depends on the mesh size $\Delta x=\Delta \hx/\hat{l}_0$ and time step $\Delta t=\Delta\hht/\hat{t}_0$ in addition to $\epsilon$, i.e., $f_h=f_h(\epsilon; \Delta x, \Delta t; \x, t)$. {Here $h$ indexes the discrete scale $h=(\Delta x,\Delta t)$}. If the numerical scheme is a good discretization of the Boltzmann equation, it can be expected that $f_h(\epsilon; \Delta x, \Delta t) \to f(\epsilon)$ as $\Delta x \to 0$ and $\Delta t\to 0$ for a fixed $\epsilon$, such that $f_h$, like $f$, is consistent with the hydrodynamic equations for small $\epsilon$. However, the computational burden may be extremely intensive if the scheme requires $\Delta x \ll \epsilon$ and/or $\Delta t \ll \epsilon$ for a valid hydrodynamic solution, such as Navier-Stokes ones. A practical kinetic scheme should be able to give a solution without resolving the kinetic scales, from which reliable numerical solutions of the hydrodynamic equations can be obtained, such as those direct Navier-Stokes solvers from discretizing the Navier-Stokes equations without imposing the mesh size and time step scales on the particle mean free path and collision time scales. Therefore, the current question is whether, and how much, $f_h$ can be consistent with the corresponding hydrodynamic equations as $\Delta x \gg \epsilon$ and $\Delta t \gg \epsilon$ (i.e., $\Delta \hx \gg \hat{\lambda}_0$ and $\Delta \hht \gg \hat{\tau}_0$).

Currently, the most popular approach to analyze the asymptotic behaviors of a kinetic scheme is to check its asymptotic preserving (AP) property, which is about to identify a consistent and stable discretization of the hydrodynamic equations in the limit of $\epsilon \to 0$ without resolving the kinetic scales \cite{ref:AP_Rev_2017,ref:AP_Rev_2012,ref:Jin1999}.
}
Historically, the study of AP properties of numerical schemes for kinetic equations can be traced back to the pioneering work of Larsen \cite{ref:Larsen_1983}, where the asymptotic behaviors of several spatial differencing schemes for the steady neutron transport equation were analyzed in the diffusion limit ($\epsilon_n\to 0$) with a fixed spatial mesh $\Delta x$, via an asymptotic expansion of the numerical scheme with respect to the neutron Knudsen number $\epsilon_n$. The methodology was systematized later by considering the diffusive asymptotic behavior of numerical schemes under different mesh scales \cite{ref:Larsen_1987}.
The approach was also employed to analyze and design numerical schemes for time-dependent transport equation with the diffusive scaling \cite{ref:Klar98} and hyperbolic systems with stiff relaxation terms \cite{ref:Jin1995,ref:Pareschi2003}. Particularly, a concrete definition of AP scheme was presented in the analysis of numerical methods for a linear hyperbolic system containing a stiff relaxation \cite{ref:Lowrie02}, focusing on the spatial discretizations.

The AP concept was also employed to analyze and develop kinetic schemes for flow problems. For instance, early in 1991, a time-splitting scheme for the Bhatnagar-Gross-Krook (BGK) equation was already proposed such that it can preserve the Euler limit \cite{ref:Coron91} under acoustic scaling, and a finite-difference scheme for the Boltzmann equation, which can capture the incompressible Navier-Stokes equations under diffusive scaling, was developed later \cite{ref:AP-NS-Klar}. Actually, a large number of kinetic schemes with AP properties have been devised in the past three decades, such as the IMEX method \cite{ref:IMEX_Hu18}, micro-macro decomposition method \cite{ref:Luc_NS_JCP08}, and penality method \cite{ref:AP-NS-FJJCP2010}. Several review papers at different periods summarized the progress on this subject \cite{ref:Review14, ref:AP_Rev_2017,ref:Luc_RGD2014}.

The literature review reveals several important points about the state of art of the AP methodology:

(i) Most studies focus on the leading order asymptotics ($\epsilon\to 0$) of kinetic methods under different scalings, such as the Euler limit under the acoustic scaling (e.g., \cite{ref:Coron91,ref:Tcher01}) or incompressible Navier-Stokes limit under the diffusive scaling (e.g., \cite{ref:AP-NS-Klar,ref:KlarJCP99}). Only a few schemes have been developed with the consideration of non-leading order (i.e., compressible Navier-Stokes) asymptotics for $0<\epsilon \ll 1$ under the acoustic scaling \cite{ref:IMEX17,ref:Luc_NS_JCP08, ref:AP-NS-FJJCP2010,ref:AP_NS_Hu2017,ref:AP_NS_Tong_JCP15}. However, it should be noted that it is implausible to use the AP framework for the analysis of non-leading order asymptotic behaviors since AP concept is defined by the limiting behavior as the scaling parameter approaches to zero \cite{ref:Jin1999}.  

(ii) In most previous studies, the time asymptotics and space asymptotics are usually considered independently even though a few fully discrete kinetic schemes were analyzed, e.g. \cite{ref:Luc_NS_JCP08,ref:AP_NS_Tong_JCP15}. Indeed, when attention is paid on the time asymptotics, the spatial gradient is assumed to be discretized exactly or the influence of spatial discretization is ignored (e.g., \cite{ref:IMEX17,ref:IMEX07,ref:IMEX_Hu18,ref:AP-NS-FJJCP2010,ref:Tcher01,ref:KlarJCP99,ref:AP_NS_Hu2017,ref:Klar_CF2006,AP-NS-Bos2017}).
In contrast, when spatial asymptotics is concerned, the semi-discrete form of kinetic equations is taken. In other words, the time integration is assumed to be precise. For example, several semi-discrete schemes for non-stationary hyperbolic systems were analyzed \cite{ref:Lowrie02,ref:Jin_Levermore1996}, where it was shown that the flux reconstruction is crtitical for space asmyptotics, as already demonstrated in \cite{ref:Larsen_1987} for the steady neutron transport equation.
For time-dependent problems, the time evolution and spatial gradients interact mutually. In practical applications, the temporal-spatial coupling feature becomes important to capture correct physics \cite{ref:Li-TS-2019,ref:Li-Du-2016,ref:Pan-Xu-Li-Li-2016}. Consequently, the temporal-spatial decoupled discretizations may lead to large numerical dissipations and thus influence the asymptotics, as argued and numerically verified in \cite{ref:ChenXu-2015} by simulating the steady cavity flow. It is noted that some kinetic schemes with temporal-spatial coupling have been developed in recent years, such as the UGKS \cite{ref:UGKS} and DUGKS \cite{ref:DUGKS13} methods. In both methods, the flux at cell interface is reconstructed based on evolution solutions of the kinetic equation itself, and thus the coupled particle transport and collision in the temporal-spatial evolution is fully incorporated into the numerics. Indeed, a number of previous studies have demonstrated the distinguished features of these methods \cite{ref:UGKS,ref:ChenXu-2015,ref:GDVM-WP}. We will further demonstrate the importance of coupled transport and collision in designing kinetic schemes in the present work.


(iii) The standard procedure of asymptotic analysis in the AP framework is to verify a numerical scheme for the kinetic equation reduces to a consistent and stable scheme of the hydrodynamic (Euler or Navier-Stokes) equations, where $\Delta x$ and $\Delta t$ are fixed and independent of the scaling parameter $\epsilon$ \cite{ref:Review14,ref:AP_Rev_2017}. Many kinetic schemes with either Euler or (in)compressible Navier-Stokes asymptotics have been analyzed following this procedure, e.g., \cite{ref:AP-NS-Klar,ref:Luc_NS_JCP08,ref:AP_NS_Tong_JCP15}. An alternative but more insightful approach for asymptotic analysis is to make use of the modified equation of the kinetic scheme, which is a standard technique in classical partial differential equations for illustrating the consistency of a numerical scheme as well as the dissipation and dispersion features \cite{ref:ME,ref:Li-2013}. For example, Jin and Levermore applied the asymptotic expansion to some semi-discrete (in space) schemes for hyperbolic systems with stiff relaxations \cite{ref:Jin_Levermore1996}, and then derived the so-called asymptotic modified equations (AME) corresponding to the limiting numerical schemes, from which the asymptotic properties were investigated by evaluating the leading order term in the truncation error under different scalings of the mesh size with respect to the relaxation parameter. Lowrie and Morel further analyzed the influence of high-order terms in the truncation error in the AME on the AP properties \cite{ref:Lowrie02}. A similar method was also employed to analyze the time asymptotics of a class of IMEX Runge-Kutta schemes at the compressible Navier-Stokes level \cite{ref:AP_NS_Hu2017}.  It is clearly shown that the modified equation approach can reveal the influence of scalings of mesh size (or time step) on the asymptotic behavior.

With the above understandings and motivated by the previous studies, in this article we attempt to present a
clear picture of the {\it asymptotic process} for fully transport-collision coupled discrete kinetic schemes at small but non-zero Knudsen number, through the modified equation approach, aiming to provide the insight into the behavior of numerical solutions at  small but non-zero mesh sizes. Specifically,
 we will introduce a new concept of {\it unified preserving} (UP) property, which is able to assess orders of asymptotics and the limiting equations underlying the scheme at small Knudsen numbers.
 { This concept is motivated by the Taylor expansion analysis which can be employed to asses the degree of approximation between two functions. Given two functions $g(x)$ and $h(x)$ with sufficient smoothness, the Taylor series of them at $x=0$ can be expressed as $g(x)=\sum_{k=0}^{\infty}{g^{(k)}x^k/k!}$ and $h(x)=\sum_{k=0}^{\infty}{h^{(k)}x^k/k!}$, where $g^{(k)}$ and $h^{(k)}$ are the  $k$-th derivative of $g$ and $h$ evaluated at the point $x=0$, respectively. Then if $g^{(k)} = h^{(k)}$ for $k=0\sim n$, but $g^{(n+1)} \neq h^{(n+1)}$, we can say $g(x)$ and $h(x)$ are consistent with each other up to order $n$. Borrowing this idea, we define the UP property of a kinetic scheme to assess the approximation degree of $f_h$ to $f$ in terms of their Chapman-Enskog expansions. Specifically, if $f_h$ can be expanded as  $f_h(\epsilon; \Delta x, \Delta t; \x, t)=\sum_{k=0}^{\infty}{\epsilon^k f_h^{(k)} (\Delta x, \Delta t; \x, t)}$, we can estimate the consistency between $f_h$ and $f$ by comparing the expansion coefficients $f^{(k)}$ and $f_h^{(k)}$ for $k=0, 1, \cdots $, such that the approximation order of $f_h$ to $f$ can be determined. Since $f$ is consistent with the hydrodynamic equations for small $\epsilon$, it is expected that a consistent numerical solution of certain hydrodynamic equations can be extracted from $f_h$, depending on the approximation order of $f_h$ to $f$. Particularly, at the leading order, $f_h$ can give a numerical solution consistent with the Euler equations, or in other words, the kinetic scheme reduces to a numerical scheme for the Euler equations as $\epsilon\to 0$. In this regard, the UP concept is consistent with AP. But at higher orders, it can distinguish the hydrodynamic behaviors beyond the Euler limit, as those by the Chapman-Enskog analysis of the Boltzmann equation. Therefore, the UP concept is a generalization of AP, in the sense that it can not only be used to asses the hydrodynamic behaviors of a kinetic scheme at the Euler limit, but also can provide its high-order hydrodynamics and estimate the order of asymptotics.
 }

The paper is organized as follows.  We summarize the kinetic equation and its asymptotic behavior at small Knudsen number in Sec. \ref{sec:KE}. The concept of unified preserving schemes is proposed and justified in Sec. \ref{sec:UP}. As two examples, the discrete unified gas-kinetic scheme (DUGKS) and a second-order IMEX scheme are analyzed in Sec. \ref{sec:DUGKS} to demonstrate their UP properties. A summary is given  in
Sec. \ref{sec:summary}. We also provide a UP analysis of a kinetic scheme with collisionless reconstruction of cell interface distribution function in \ref{sec:appen} to show its dynamic difference from the DUGKS in the hydrodynamic regime.

\section{Kinetic equation and asymptotic behaviors}
\label{sec:KE}

In this section we first give a brief introduction of kinetic models and the asymptotic analysis method which will be used later in this work. Although these are standard and can be found in books, e.g., \cite{ref:ChapmanBook,ref:CercignaniBook}, essential details presented here will make the later analysis much clearer. For simplicity we consider the BGK equation for a monatomic gas without external force (the same analysis should apply for general cases),
\begin{equation}
\label{eq:BGK}
\partial_{\hht}{\hf}+\hxi\cdot\hat{\nabla} \hf = \hat{Q},
\end{equation}
where $\partial_\hht=\partial/\partial{\hht}$ and $\hat{\nabla}=\partial/\partial{\hx}$, $f(\hx, \hxi, \hht)$ is the distribution function at time $\hht$ and position $\hx$ for particles moving with velocity $\hxi$, and the collision operator is
\begin{equation}
\hQ =-\dfrac{1}{\htau}\left(\hf-\hf^{(eq)}\right),
\end{equation}
with $\htau$ being the relaxation time which depends on the pressure and viscosity coefficient, and $\hf^{(eq)}$ is the local Maxwellian equilibrium defined by the hydrodynamic variables $\hat{\bm{W}}=(\hrho, \hu, \hT)^T$,
\begin{equation}
\hf^{(eq)}(\hat{\bm{W}})=\dfrac{\hrho}{(2\pi R\hT)^{D/2}}\exp\left(-\dfrac{|\hxi-\hu|^2}{2R\hT}\right),
\end{equation}
where $D$ is the spatial dimension, $R$ is gas constant, while $\hrho$, $\hu$, and $\hT$ are the gas density, velocity and temperature, respectively,
\begin{equation}
\label{eq:W}
\hat{\bm{W}}=\left(
\begin{array}{c}
\hrho \\
\hu \\
\hT
\end{array}
\right)
=
\left(
\begin{array}{c}
\int{\hf\ \d\hxi} \\
\dfrac{1}{\hrho}\int{\hxi \hf\ \d\hxi} \\
\dfrac{1}{D\hrho R}\int{|\hxi-\hu|^2 \hf\ \d\hxi}
\end{array}
\right).
\end{equation}
Note that the BGK operator is collisional invariant, i.e.,
\begin{equation}
\label{eq:ConservQ}
\int{\hQ\ \d\hxi}=0,\quad \int{\hxi \hQ\ \d\hxi} = \bm{0},\quad \int{|\hxi|^2 \hQ\ \d\hxi} = 0.
\end{equation}

The asymptotic behavior of the kinetic equation \eqref{eq:BGK} for small $\epsilon$ can be analyzed by the Chapman-Enskog expansion method in terms of the small parameter $\epsilon$, proportional to the non-dimensional Knudsen number introduced below, in order to relate with hydrodynamics \cite{ref:ChapmanBook}. For this purpose we first rewrite the BGK equation \eqref{eq:BGK} in a non-dimensional form by introducing the following dimensionless variables,
\begin{equation}
\label{eq:dimensionlessParameter}
\rho=\dfrac{\hrho}{\hrho_0}, \quad \u=\dfrac{\hu}{\hc_0}, \quad T=\dfrac{\hT}{\hT_0}, \quad f=\dfrac{\hf}{\hrho_0/\hc_0^D},\quad \x=\dfrac{\hx}{\hat{l}_0}, \quad t=\dfrac{\hht}{\hht_0}, \quad \bxi=\dfrac{\hxi}{\hc_0}, \quad \tau=\dfrac{\htau}{\htau_0},
\end{equation}
 where $\hrho_0$, $\hT_0$ and $\hc_0=\sqrt{2R\hT_0}$ are the reference density, temperature, and molecular velocity, respectively, while $\hat{l}_0$, $\hht_0=\hat{l}_0/\hc_0$, and $\htau_0$ are the reference length, time, and mean free time, respectively. It is noted that the relaxation time and mean free path can be related to the dynamic viscosity $\hat{\mu}$ and pressure $\hat{p}=\hrho R \hT$ \cite{ref:CercignaniBook}, namely, $\htau=\hat{\mu}/\hat{p}$ and $\hat{\lambda}=\htau\sqrt{\pi R\hT/2}$, therefore the parameter $\epsilon=\htau_0/\hht_0 = \hat{\lambda}_0 /\hat{l}_0$ (here $\hat{\lambda}_0=\hc_0 \htau_0$) is proportional to the Knudsen number $\mbox{Kn}$ with the same order, which measures the ratio between the kinetic scale $(\hat{\tau}_0, \hat{\lambda}_0)$ and the hydrodynamic scale $(\hat{t}_0, \hat{l}_0)$. The dimensionless form of Eq. \eqref{eq:BGK} can then be expressed as
\begin{equation}
\label{eq:nonDimensional-BGK}
\partial_{t}{f}+\bxi\cdot\nabla f =-\dfrac{1}{\epsilon}Q =-\dfrac{1}{\epsilon\tau}\left(f-f^{(eq)}(\rho,\u,T)\right),
\end{equation}
where
\begin{equation}
f^{(eq)}=\dfrac{\rho}{(2\pi R_0 T)^{D/2}}\exp\left(-\dfrac{|\bxi-\u|^2}{2R_0T}\right), \qquad (R_0=1/2).
\end{equation}

{
In addition to the physical variables  $\x$, $t$, and $\bxi$, the solution of Eq. \eqref{eq:nonDimensional-BGK} also depends on the parameter $\epsilon$, i.e., $f=f(\epsilon;\x,\bxi, t)$.
}
In the Chapman-Enskog analysis, it is assumed that the distribution function depends on space and time only through a functional dependence on the hydrodynamic variables, i.e., $f(\x,\bxi,t)=f(\bxi, \bm{W}(\x,t), \nabla\bm{W}(\x,t), \nabla\nabla\bm{W}(\x,t),\cdots)$. Under such assumption, the distribution function can be expressed as a series expansion in powers of $\epsilon$ \cite{ref:CercignaniBook},
\begin{equation}
\label{eq:fCE}
f=f^{(0)} + \epsilon f^{(1)} + \epsilon^2 f^{(2)} + \cdots,
\end{equation}
where the expansion coefficients $f^{(k)}$ depend on the {hydrodynamic variables $\bm{W}$} and their gradients, with the assumption that $f^{(eq)}=O(1)$ and $f^{(k)})=O(1)$ for $k\ge 0$. Correspondingly, the time derivative is also expanded formally as a series of $\epsilon$,
\begin{equation}
\label{eq:tCE}
\partial_{t}=\partial_{t_0} + \epsilon \partial_{t_1} + \epsilon^2 \partial_{t_2} + \cdots,
\end{equation}
where $\partial_{t_k}$ denotes the contribution to $\partial_t$ from the spatial gradients of the hydrodynamic variables \cite{ref:CercignaniBook,ref:ChapmanBook}. Specifically, the perturbation expansion \eqref{eq:fCE} generates similar expansions of the pressure tensor and heat flux in the hydrodynamic balance equations, and $\partial_{t_k}$ is defined to balance these terms at different orders of $\epsilon$. In general, $\partial_{t_k}$ is related to the $(k+1)$-order spatial gradients of the hydrodynamic variables.

The expansion coefficients $f^{(k)}$ can be found by substituting the above expansions into Eq. \eqref{eq:nonDimensional-BGK} and multiplying $\epsilon$ on both sides, which gives that
\begin{subequations}
\label{eq:CE-BGK}
\begin{align}
\epsilon^0:&\qquad f^{(0)} = f^{(eq)},\\
\epsilon^1:&\qquad  D_0 f^{(0)} =   -\dfrac{1}{\tau} f^{(1)},\label{eq:CE-1}\\
\epsilon^2:&\qquad \partial_{t_1}{f^{(0)}}+D_0 f^{(1)} = -\dfrac{1}{\tau} f^{(2)}, \label{eq:CE-2} \\
\cdots & \cdots \\
\epsilon^{k}:&\qquad \sum_{j=1}^{k-1}\partial_{t_j}{f^{(k-j-1)}}+D_0 f^{(k-1)} =  -\dfrac{1}{\tau} f^{(k)}, \label{eq:CE-k}
\end{align}
\end{subequations}
where $D_0=\partial_{t_0}+\bxi\cdot\nabla$. {In this work, we call the equations at different orders of $\epsilon$ as {\it balance equations}, which connect the coefficient $f^{(k)}$ to low order ones from $f^{(0)}$ to $f^{(k-1)}$}. Noting that with the conservation property \eqref{eq:ConservQ} of the collision operator, we have
\begin{equation}
\int{\bm{\psi} f^{(k)}\ \d\bxi}=0, \quad k>0.
\end{equation}
where $\bm{\psi}=(1, \bxi, |\bxi|^2/2)$ are the collision invariants. Then the hydrodynamic equations can be derived with different approximation orders of $\epsilon$. For instance, taking the conservative moments of Eq. \eqref{eq:CE-1} leads to
\begin{equation}
\partial_{t_0} \int{\bm{\psi} f^{(0)}\ \d\bxi} + \nabla\cdot \int{\bxi \bm{\psi} f^{(0)}\ \d\bxi} =  0.
\end{equation}
Since $f^{(0)}=f^{(eq)}$ as given in Eq. \eqref{eq:CE-1}, the equations can be rewritten as
\begin{subequations}
\label{eq:Euler}
\begin{equation}
\partial_{t_0}\rho + \nabla\cdot (\rho\u) =  0,
\end{equation}
\begin{equation}
\partial_{t_0}(\rho\u) + \nabla\cdot (\rho\u\u + p\bm{I}) =  0,
\end{equation}
\begin{equation}
\partial_{t_0}(\rho E) + \nabla\cdot [(\rho E + p)\u] =  0,
\end{equation}
\end{subequations}
where $\bm{I}$ is the second-order unit tensor, $p=\rho R_0 T$ is the dimensionless pressure and $E=c_vT+\tfrac{1}{2}|\u|^2$ is the dimensionless total energy with $c_v=DR_0/2$, which are just the Euler equations when we take the first-order approximation, i.e., $\partial_t=\partial_{t_0}$.

If we further take the conservative moments of Eq. \eqref{eq:CE-2}, we can obtain
\begin{subequations}
\label{eq:H-t1}
\begin{equation}
\partial_{t_1}\rho =  0,
\end{equation}
\begin{equation}
\partial_{t_1}(\rho\u) + \nabla\cdot \bm{P}^{(1)} =  0,
\end{equation}
\begin{equation}
\partial_{t_1}(\rho E) + \nabla\cdot \bm{Q}^{(1)} =  0,
\end{equation}
\end{subequations}
where $\bm{P}^{(1)}=\int{\bxi\bxi f^{(1)}\ \d\bxi}$ and $\bm{Q}^{(1)}=\tfrac{1}{2}\int{|\bxi|^2\bxi f^{(1)}\ \d\bxi}$. From Eq. \eqref{eq:CE-1}, we can evaluate the second-order tensor $\bm{P}^{(1)}$ and the vector $\bm{Q}^{(1)}$ explicitly,
\begin{subequations}
\begin{equation}
-\dfrac{1}{\tau}P_{\alpha\beta}^{(1)} = \partial_{t_0}\left(\rho u_\alpha u_\beta + p\delta_{\alpha\beta}\right)+\nabla_\gamma{\Gamma^{(0)}_{\alpha\beta\gamma}} ,
\end{equation}
\begin{equation}
-\dfrac{1}{\tau}Q_{\alpha}^{(1)} = \partial_{t_0}\left[(p+\rho E) u_\alpha\right]+\nabla_\beta{\Theta^{(0)}_{\alpha\beta}} ,
\end{equation}
\end{subequations}
where  $\bm{\Gamma}^{(0)}=\int{\bxi \bxi\bxi f^{(0)}\ \d\bxi}$ and $\bm{\Theta}^{(0)}=\tfrac{1}{2}\int{|\bxi|^2\bxi \bxi f^{(0)}\ \d\bxi}$. The time derivatives of $\partial_{t_0}$ can be evaluated from Eq.\eqref{eq:Euler}, and after some standard algebra we can obtain that
\begin{subequations}
\begin{equation}
\label{eq:PQ}
P_{\alpha\beta}^{(1)} = -\sigma_{\alpha\beta}\equiv -\mu\left[\partial_\alpha u_\beta + \partial_\beta u_\alpha - \dfrac{2}{D}(\nabla\cdot\u)\delta_{\alpha\beta}\right] ,
\end{equation}
\begin{equation}
Q_{\alpha}^{(1)} = -\kappa\partial_\alpha T - u_\beta\sigma_{\alpha\beta},
\end{equation}
\end{subequations}
where $\mu=\tau p$ and $\kappa=\tau p (D+2)R_0/2$. Then with Eqs. \eqref{eq:Euler} and \eqref{eq:H-t1}, we  obtain the hydrodynamic equations up to $\epsilon$, i.e., $\partial_t=\partial_{t_0}+\epsilon\partial_{t_1}$,
\begin{subequations}
\label{eq:NS}
\begin{equation}
\partial_{t}\rho + \nabla\cdot (\rho\u) =  0,
\end{equation}
\begin{equation}
\partial_{t}(\rho\u) + \nabla\cdot (\rho\u\u + p\bm{I}) = \nabla\cdot(\epsilon \bm{\sigma}),
\end{equation}
\begin{equation}
\partial_{t}(\rho E) + \nabla\cdot [(\rho E + p)\u] =  \nabla\cdot(\epsilon\kappa T)+\nabla\cdot(\epsilon\bm{\sigma}\cdot\u),
\end{equation}
\end{subequations}
which are exactly the Navier-Stokes equations with the unit Prandtl number. The high-order hydrodynamic equations, such as the Burnett and super-Burnett equations, can be derived by invoking the corresponding high-order kinetic equations in the successive system \eqref{eq:CE-BGK}.

Although the Chapman-Enskog expansion method is criticized from different viewpoints \cite{ref:Gorban2014}, it is still a useful approach to analyze the asymptotic behavior of kinetic models at small $\epsilon$. Particularly, the Euler and Navier-Stokes equations derived from the kinetic model can precisely describe the corresponding hydrodynamic behaviors of the kinetic model.
{The higher order hydrodynamic equations obtained from the Chapman-Enskog expansions usually suffer from physical and numerical difficulties. In order to overcome this defect, some modified versions of the Chapman-Enskog expansion have been proposed, such as the regularization one \cite{ref:RegCE_Rose1989}, the renormalization one \cite{ref:NormCE_Slemrod1989}, and the exact summation one \cite{ref:SumCE_Gorban1982}. With these reformulated Chapman-Enskog expansions, hydrodynamic equations with desired physics can be derived from the kinetic equation at higher orders, and it is now recognized that ``{\it\textbf{the problem with the Chapman-Enskog expansion is not the expansion itself but its truncation}}" \cite{ref:NormCE_Slemrod2012}. Therefore, the Chapman-Enskog expansion can still serve as a basis for analyzing the asymptotic behaviors of kinetic models as well as kinetic schemes. For simplicity, in the present work we will adopt the original Chapman-Enskog expansion method to illustrate the basic idea, but the analysis procedure presented here can also be conducted based on the modified versions.
Furthermore, we note that high-order Chapman-Enskog analysis beyond the Navier-Stokes order is also helpful for designing kinetic schemes. For example, some lattice Boltzmann models incorporating the Burnett effect have been developed for two-phase flows with better performance \cite{ref:LB_Wagner2006,ref:LB_HuangRZ2016,ref:LB_ZhengL2017}, and a gas-kinetic scheme at the Burnett level has been developed for gas flows in slip regime \cite{ref:BGK_Burnett}.
}

\section{Definition of Unified Preserving (UP) Property}
\label{sec:UP}
A numerical scheme for the kinetic equation \eqref{eq:nonDimensional-BGK}, denoted by $P_h^{\epsilon}$, gives an approximate distribution function $f_h$,  depending on the cell size $\Delta x$ and time step $\Delta t$, i.e., $f_h=f_h(\epsilon; \x, \bxi, t;\Delta x, \Delta t)$.
Therefore, it is expected that the asymptotic behavior of $P_h^{\epsilon}$ depends not only on the kinetic and hydrodynamics scales, but also on the numerical scale $h=(\Delta x, \Delta t)$ (or $\Delta \hx$ and $\Delta \hht$ in dimensional form). It is noted that we concentrate here on the time and space discretizations while the velocity discretization is not considered in this work. Furthermore, we assume that the cell size and time step are adequate to resolve the flow physics at the hydrodynamic scale, similar to many Navier-Stokes solvers.

In order to assess the asymptotic property of $P_h^{\epsilon}$ at small $\epsilon$ on the hydrodynamic scale, we can again apply the Chapman-Enskog expansion to $f_h$,
\begin{equation}
\label{eq:fhCE}
f_h=f_h^{(0)} + \epsilon f_h^{(1)} + \epsilon^2 f_h^{(2)} + \cdots.
\end{equation}
Then we  compare the expansion coefficients $f_h^{(k)}$ with those of the original distribution function $f(\x, \bxi, t)$ given by Eq. \eqref{eq:fCE}. The comparison leads to the definition of {\it unified preserving} property of the kinetic scheme $P_h^{\epsilon}$.

\begin{definition}
A consistent numerical scheme of Eq. \eqref{eq:nonDimensional-BGK}, $P_h^{\epsilon}$, is called an $n$-th order {\it unified preserving} (UP) scheme
provided that
\begin{enumerate}
         \item[(i)] it is uniformly stable, i.e., the scheme is stable regardless of $\epsilon$;\\
         {\item[(ii)] it is a legitimate discretization of the collisionless kinetic equation as $\epsilon \to \infty$; }\\
         \item[(iii)] for $\epsilon \ll 1$,  there exist two parameters { $\alpha_0\in [0, 1)$ and $\beta_0 \in [0, 1)$, such that as $\Delta t=O(\epsilon^\alpha)$ and $\Delta x=O(\epsilon^\beta)$ with $\alpha_0<\alpha<1$ and $\beta_0<\beta<1$}, the expansion coefficients and the associated discrete balance equations satisfy
\begin{equation}
\label{eq:UP}
f_h^{(0)}=f_h^{(eq)},\quad
\sum_{j=1}^{k-1}\partial^{(j)}_{t}{f_h^{(k-j-1)}}+D_0 f_h^{(k-1)} = - \dfrac{1}{\tau}f_h^{(k)}\quad  (1\le k\le n),
\end{equation}
  but $f_h^{(n+1)}$ depends on $\Delta x$ or $\Delta t$ explicitly.
\end{enumerate}
\end{definition}

In the above definition, $f_h^{(eq)}=f^{(eq)}(\bm{W}_h)$ is the Maxwellian distribution function dependent on the numerical hydrodynamic variables $\bm{W}_h=(\rho_h, \u_h, T_h)^T$, which is defined as in Eq. \eqref{eq:W} with $f$ being replaced by $f_h$. 

The first (stable) condition suggests that the time step is not limited by the relaxation time in terms of numerical stability, which requires a non-explicit treatment of the (stiff) collision term.
{
The second condition suggests that the scheme can capture free transport of particles in the free molecular regime, which is essential for a multiscale kinetic scheme.
}

The third (asymptotic) condition indicates that the scheme can capture the hydrodynamic behaviors to some degree without resolving the kinetic scales. To see this more clearly, we introduce two parameters,
\begin{equation}
\delta_t\equiv\frac{\Delta t}{\epsilon}=\frac{\Delta \hat{t}}{\hat{\tau}_0}, \quad \delta_x\equiv\frac{\Delta x}{\epsilon}=\frac{\Delta \hat{x}}{\hat{\lambda}_0},
\end{equation}
which measure the numerical temporal and spatial resolutions of the kinetic scales. For sufficient small $\epsilon$, it can be seen that
$\delta_t=O(\epsilon^{\alpha-1}) \gg 1$ and $\delta_x=O(\epsilon^{\beta-1}) \gg 1$, meaning that the numerical scale is much larger than the kinetic scale.
Furthermore, the condition $f_h^{(0)}=f_h^{(eq)}$ suggests that the discrete collision operator is also conservative, and therefore
\begin{equation}
\int{\bm{\psi} f_h^{(k)}\ \d\bxi}=0, \quad k>0.
\end{equation}
Then it can be seen that the second equation in \eqref{eq:UP} has the same property as Eq. \eqref{eq:CE-k}, and their moment equations are also identical for $1\le k\le n$. For instance, for a first-order ($n=1$) UP scheme, we have
\begin{equation}
\label{eq:hCE-1}
D_0 f_h^{(0)}=-\dfrac{1}{\tau}f_h^{(1)}.
\end{equation}
Taking conservative moments of the above equation leads to
\begin{subequations}
\label{eq:hEuler}
\begin{equation}
\partial_{t_0}\rho_h + \nabla\cdot (\rho_h\u_h) =  0,
\end{equation}
\begin{equation}
\partial_{t_0}(\rho_h\u_h) + \nabla\cdot (\rho_h\u_h\u_h + p_h\bm{I}) =  0,
\end{equation}
\begin{equation}
\partial_{t_0}(\rho_h E_h) + \nabla\cdot [(\rho_h E_h + p_h)\u_h] =  0,
\end{equation}
\end{subequations}
where $p_h=\rho_h R_0 T_h$ and $E_h=c_vT_h+\tfrac{1}{2}|\u_h|^2$. It can be seen that Eqs. \eqref{eq:hEuler} are the same as those of the original BGK equation at the Euler level (see Eqs. \eqref{eq:Euler}), which means that the numerical hydrodynamic quantities $\bm{W}_h$ are the solutions of the Euler equations. In other words, a first-order UP kinetic scheme can reproduce the Euler equations exactly with a coarse numerical resolution without resolving the kinetic scale, and this fact also indicates that the scheme has the AP properties. However,
since the equation  of a first-order UP scheme for $f_h^{(2)}$ is incomparable with that of $f^{(2)}$, the balance moment equations at $\epsilon^2$ are different, such that the Navier-Stokes equations cannot be recovered from the scheme. Actually, in order to capture the hydrodynamic behaviors at the Navier-Stokes level, a second-order UP scheme is required. In this case, $f_h^{(1)}$ satisfies Eq. \eqref{eq:hCE-1}, and $f_h^{(2)}$ satisfies
\begin{equation}
\label{eq:hCE-2}
\partial_{t_1}{f_h^{(0)}}+D_0 f_h^{(1)} = -\dfrac{1}{\tau} f_h^{(2)},
\end{equation}
which is the same as Eq. \eqref{eq:CE-1}, and the conservative moment equations take the same form as Eq. \eqref{eq:H-t1} together with Eq. \eqref{eq:PQ}, with $\bm{W}$ being replaced by $\bm{W}_h$. As such, the numerical quantities $\bm{W}_h$ from the second-order UP scheme satisfy the Navier-Stokes equations. To proceed further, a third-order UP scheme can capture the hydrodynamic behaviors at the Burnett level, and a fourth-order one will go to the super-Burnett level. Generally, for an $n$-th UP scheme, the expansion coefficients of $f_h^{(k)}$ satisfy the same equations of $f^{(k)}$ for $0\le k \le n$, such that the hydrodynamic quantities $\bm{W}_h$ also satisfy the corresponding $n$-th order hydrodynamic equations obtained from the Chapman-Enskog analysis of the original kinetic equation. Therefore, the concept of UP can be used to distinguish the asymptotic limiting equations of different kinetic schemes for small Knudsen numbers with different orders.

\begin{remark}
It should be noted that a UP scheme is a consistent discretization of the kinetic equation \eqref{eq:BGK}, and $f_h$ is an approximation solution of the kinetic equation rather than the solution of the limiting hydrodynamic equations. Therefore, for relative large $\epsilon$, $f_h$ can still be a good approximation to the solution of the kinetic equation instead of the high-order hydrodynamic equations such as the Burnett or super-Burnett ones.
\end{remark}

\begin{remark}
It is noted that $\Delta t = \delta_t\epsilon=\Delta \hat{t}/\hat{t}_0$ and $\Delta x = \delta_x\epsilon=\Delta \hat{x}/\hat{l}_0$, which represent the numerical time and mesh resolutions for the hydrodynamic scale. Generally, the numerical resolutions must resolve the hydrodynamic scale flow sctructure for hydrodynamic problems, i.e., $\Delta t <1 $ and  $\Delta x <1$, suggesting that { $\alpha \ge 0$} and { $\beta \ge 0$}. {Particularly, as $\alpha=0$ and/or $\beta=0$, we have $\delta_t = O(\epsilon^{-1})$ and/or $\delta_x = O(\epsilon^{-1})$, meaning that the kinetic time and/or length scales are unresolved, and this is referred to as ``unresolved regime" in \cite{ref:Lowrie02} or ``thick regime" in \cite{ref:Larsen_1987}. However, this does not mean that the hydrodynamic scales are not resolved. Generally, if the flow exhibits certain hydrodynamic structure, the numerical resolutions should resolve the typical hydrodynamic scales in order to capture the flow structure, such as the hydrodynamic boundary layer.}
\end{remark}

\begin{figure}[!htb]
\centering
\includegraphics[width=0.45\textwidth]{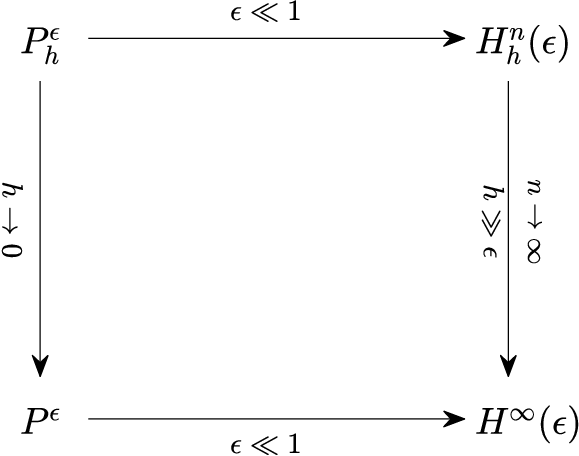}
\caption{Schematic diagram of a UP kinetic scheme of order $n$. Here $P^\epsilon$ is the kinetic equation, $P_h^\epsilon$ is a consistent kinetic scheme with discrete resolution $h$, $H^\infty(\epsilon)$ is the set of the Chapman-Enskog expansion coefficients of the solution of $P^\epsilon$, and $H_h^{k}(\epsilon)$ is the set of the expansion coefficients of the solution of $P_h^\epsilon$. For an $n$-th order UP scheme, $f_h^{(k)}=f^{(k)}$ for $0\le k \le n$.}
\label{fig:upFig}
\end{figure}

The UP property of $P_h^\epsilon$ can be analyzed by studying the asymptotics of its modified (or equivalent) equation. We consider the one-dimensional case without loss of generality, then the modified equation of a numerical scheme for \eqref{eq:nonDimensional-BGK} can generally be expressed as
\begin{equation}
\label{eq:nModifiedEq}
\partial_t{f_h}+\xi\partial_x f_h +O(\Delta t^s,\Delta x^q){\cal L} f_h
 =  \dfrac{1}{\epsilon}Q(f_h)+ \dfrac{1}{\epsilon} O(\Delta t^r) {\cal L} Q_h,
\end{equation}
where $s$, $q$, and $r$ are the orders of the leading terms of the truncation errors for the temporal integration, spatial gradient discretization, and collision integration, respectively. Note that the  time derivatives of $f_h$ with orders higher than one can be eliminated by using the equation itself.
Then substituting the Chapman-Enskog expansion of $f_h$ given by Eq. \eqref{eq:fhCE} and the time expansion \eqref{eq:tCE} into Eq. \eqref{eq:nModifiedEq}, one can obtain the set of successive equations of $f_h^{(k)}$ in terms of the order of $\epsilon$. Note that we can substitute them successively in terms of $\partial_x^q$ terms, which is just the same strategy as adopted in the modified equation approach \cite{ref:ME,ref:Li-2013}.

{
Substituting the expansions \eqref{eq:tCE} and \eqref{eq:fhCE} into the modified equation \eqref{eq:nModifiedEq}, together with the numerical scales $\Delta t =O(\epsilon^\alpha)$ and $\Delta x =O(\epsilon^\beta)$, we can obtain the expansion coefficients $f_h^{(k)}$ and the corresponding balance equations at different orders of $\epsilon$. The UP order can then be determined by comparing $f_h^{(k)}$ with $f^{(k)}$. In general, for an $n$-th order UP scheme, it is required that $\Delta t^s=o(\epsilon^{n-1})$, $\Delta t^r=o(\epsilon^{n-1})$, and $\Delta x^q=o(\epsilon^{n-1})$, since in this case the terms associated with $\Delta x$ and $\Delta t$ do not appear in the balance equations for $f_h^{(k)}$ for $0\le k\le n$.
}

It is clear that the UP order of $P_h^{\epsilon}$ depends on not only the numerical scale $h=(\Delta x, \Delta t)$, but also the discretization accuracy $(s, q, r)$. Physically, this indicates the capability of the scheme in resolving the flow physics depends on not only the numerical resolution (time step and grid spacing) or the grid Knudsen number ($1/\delta_x$), as argued in \cite{ref:Xu-POF2017}, but also the accuracy of the scheme. In other words, for a given mesh size and time step the numerical flow physics can be different from different kinetic schemes with different accuracy.

The UP property of a kinetic scheme is illustrated in Fig. \ref{fig:upFig}. The original kinetic equation with the small parameter $\epsilon$ is represented by $P^{\epsilon}$, and $P_h^{\epsilon}$ is a consistent discretization of $P^{\epsilon}$ with numerical scale $h$, namely, $P_h^{\epsilon}\to P^{\epsilon}$ as $h\to 0$ (represented by the left downward arrow). $H^{\infty}(\epsilon)=\{f^{(k)}| 0\le k <\infty \}$ is the set of the Chapman-Enskog expansion coefficients of $f$ defined by Eq. \eqref{eq:CE-BGK}, which is determined from the kinetic equation $P^{\epsilon}$ at small $\epsilon$ (represented by the bottom rightward arrow). $H_h^{n}(\epsilon)=\{f_h^{(k)}| 0\le k \le n \}$ is the set of the Chapman-Enskog expansion coefficients of $f_h$ defined by Eq. \eqref{eq:UP}, which is determined from the kinetic scheme $P_h^{\epsilon}$ at small $\epsilon$ (represented by the top rightward arrow). An $n$-th order UP scheme means that $H_h^{n}(\epsilon)$ is a subset of $H^{\infty}(\epsilon)$ as $h\gg \epsilon$ and approaches to $H^{\infty}(\epsilon)$ closer with increasing $n$, which is represented by the right downward arrow.
Note that $H^{\infty}(\epsilon)$ and $H_h^{n}(\epsilon)$ also determine uniquely the moment equations of the kinetic equation $P^{\epsilon}$ and the kinetic scheme $P_h^{n}(\epsilon)$ at different orders of $\epsilon$, respectively, and therefore the underlying moment equations of $P_h^{\epsilon}$ at orders of $\epsilon$ from 0 to $n$ are the same as those of $P^{\epsilon}$.

\begin{figure}[!htb]
	\centering
	\includegraphics[width=0.6\textwidth]{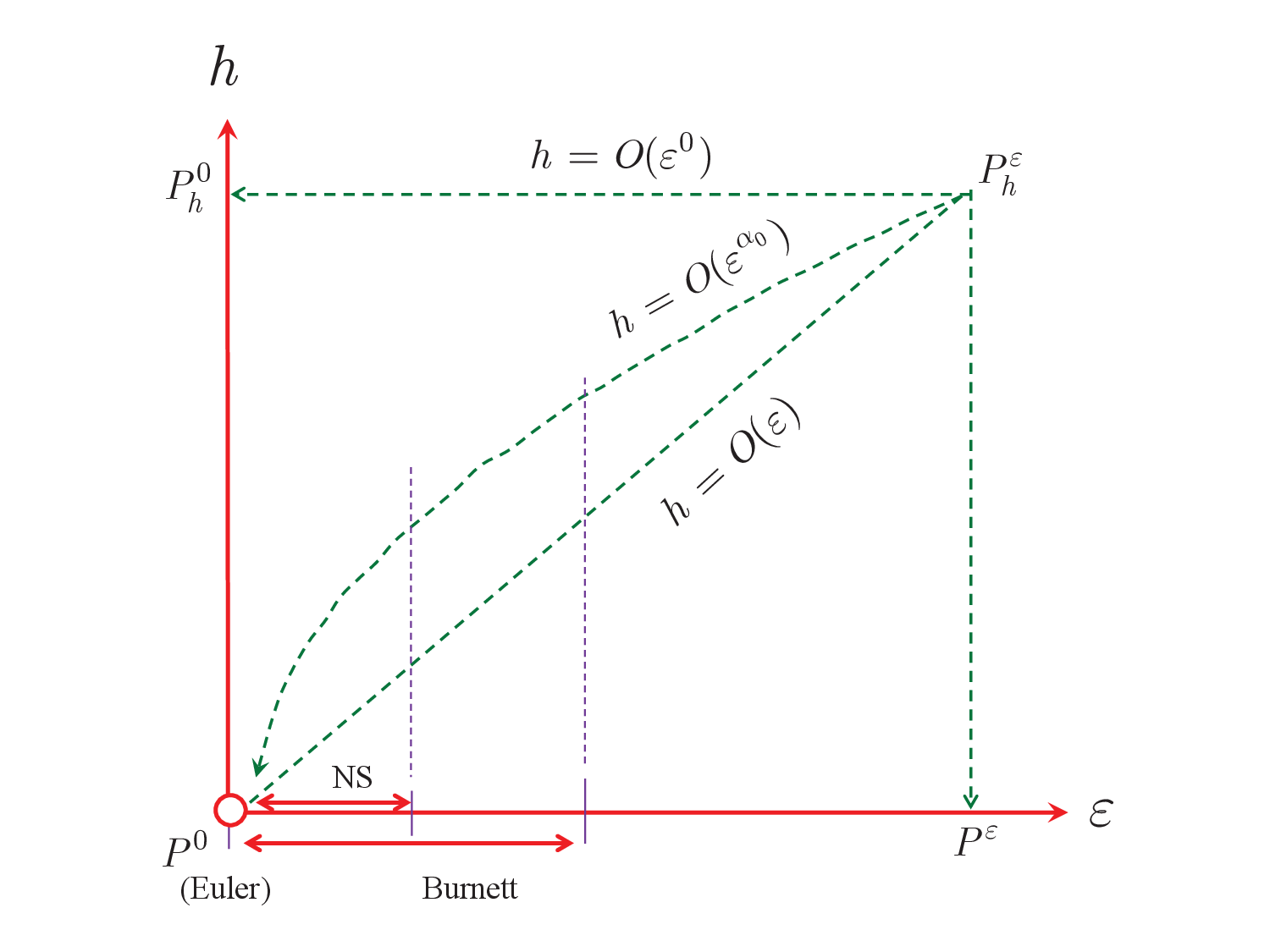}
	\caption{Schematic diagram of the asymptotic path to the limiting hydrodynamic regimes. The region below the line $h=O(\epsilon)$ suggests resolved kinetic scale. The line $h=O(\epsilon^{\alpha_0})$ represents the upper limit for a UP scheme. The region between the two lines $h=O(\epsilon)$ and  $h=O(\epsilon^{\alpha_0})$ represents the parameter space for a UP scheme.}
	\label{fig:upPath}
\end{figure}

In addition to the function on the assessment of asymptotic orders, the UP concept also provides a picture how a kinetic scheme approaches to the asymptotic limit, as demonstrated in Fig. \ref{fig:upPath}.  In the classical AP concept, a kinetic scheme $P_h^{\epsilon}$ approaches to its discrete asymptotic limit $P_h^0$ first following the line of $h=O(\epsilon^0)$ (i.e., $\epsilon\to 0$ with fixed $\Delta x$ and/or $\Delta t$), and then approaches to the analytical limit $P^0$ following the $h$-axis ($h\to 0$). This suggests that an asymptotic scheme may not pass the Navier-Stokes region, as demonstrated in \cite{ref:ChenXu-2015} where an AP-IMEX scheme was shown to fail to capture the Navier-Stokes solution of a steady flow.
On the other hand, a UP scheme of order higher than 2 can approach to the hydrodynamic limit passing through the Navier-Stokes regime following a curve below the line $h=O(\epsilon^{\alpha_0})$.

{
	 We remark that the modified equation approach was also employed in some previous AP analysis, as noted in item (iii) of the introduction section. Particularly, Lowrie and Morel studied the asymptotics behavior of a certain space semi-discrete numerical scheme for hyperbolic systems with a stiff relaxation term \cite{ref:Lowrie02}, by applying the asymptotic expansion to the modified equation of the scheme. A so-called {\it asymptotic modified equation} (AME) is then derived from the multiscale expansion, which can be viewed as the modified equation for the limiting asymptotic equation. By examining the truncation error that may depend on the scaling parameter $\epsilon$ and mesh size $\Delta x$, the leading-order asymptotic behavior of the numerical scheme can be figured out. It was shown that the dimensionless parameter characteristizing the near-equilibrium limit should be fixed to discrimate AP and non-AP schemes, and the high-order error terms in $\Delta x$ dominate the AP property.
	The UP analysis shares some similarities with the above AME approach, in the sense that both analyze the asymptotic behavior of a numerical scheme by applying the asymptotic expansion to the modified equation. However, some discrepancies between the two approaches are also clear. Firstly, the aim of the AME approach is to distinguish whether or not a scheme is AP for the leading-order solution in the limit of $\epsilon \to 0$, while the UP aims to figure out the asymptotic behaviors at different orders of $\epsilon$, which presents more detailed AP information quantitatively. Secondly, in the AME approach the asymptotic modified equation corresponding to the numerical scheme for the limiting asymptotic equation is derived explicitly, from which the influence of the mesh size on the AP can be analyzed. On the other hand, in the UP analysis only the modified equation for the original numerical scheme \eqref{eq:nModifiedEq} is presented, which is used only to derive the expansion coefficients (and the balance equation at each order), while the modified equation for the corresponding asymptotic numerical scheme does not appear. Finally, unlike the UP analysis, in the AME approach the numerical solution ($f_h$) is not distinguished from the solution of the original continuous equation ($f$), which may lead to some confusions between their asymptotic behaviors.
}

{
	It is noted that the discrete Chapman-Enskog expansion was also adopted in some previous studies on time asymptotics of certain semi-discrete kinetic schemes, e.g., \cite{ref:IMEX17,ref:Luc_NS_JCP08,ref:AP-NS-FJJCP2010,ref:AP_NS_Hu2017,ref:AP_NS_Tong_JCP15}. However, the analysis in these studies focused on time asymptotics where the spatial gradient of the scheme was assumed to be exact. As will be shown later in the UP analysis of several example kinetic schemes, the spatial discretization plays a key role for their aysmptotic behaviors, which means the sole time asymptotic analysis is insufficient for evaluating the overall behaviors at different hydrodynamic levels. This is reasonable since the temporal-spatial coupling reflects the real physics underlying the kinetic equation. Furthermore, in these studies the Chapman-Enskog approximation of the discrete distribution function was substituted into the numerical scheme. The analysis shows that the moments of the numerical scheme are just certain discretization formulations of the corresponding hydrodynamic equations. On the other hand, the Chapman-Enskog expansion of the numerical solution in the UP analysis is applied to the modified equation, such that the hydrodynamic equations underlying the schemes with different numerical resolutions can be revealed, which are unavilable in the former approach.
}

{
Finally, it is interesting to compare the present UP analysis with the multiscale expansion (ME) method for lattice Boltzmann equation (LBE) (e.g. \cite{ref:Junk-Yong2002,ref:Junk2005}), which is widely used for analyzing the hydrodynamic behaviors of LBE. As noted above, in the UP analysis three scale levels are involved, i.e., kinetic scales $(\hat\lambda_0, \hat\tau_0)$, numerical scales $(\Delta\hx,\Delta\hht)$, and hydrodynamics scales $(\hat{l}_0, \hht_0)$. The hydrodynamic behaviors of a kinetic scheme depend not only on the kinetic and hydrodynamic scales, measured by the Knudsen number $\epsilon$ as the scaling parameter, but also on the numerical scales.
On the other hand, in the ME analysis \cite{ref:Junk2005}, the LBE is just viewed as a kinetic model similar to the continuous Boltzmann equation, in which the lattice time step and lattice spacing are taken as the collision time and mean-free path. In other words, in LBE the numerical scales are the same as kinetic scales. This is also evident from the definition of the scaling parameter ($\epsilon=\Delta x \equiv \Delta \hx/\hat{l}_0$) used in the expansion, which is clearly different from the Knudsen number used in the UP analysis. The aim of the multiscale analysis of LBE is to show its hydrodynamic behaviors, just as the original Chapman-Enskog or asymptotic analysis for the continuous Boltzmann equation. Therefore, in the multiscale analysis used for LBE, the physical kinetic scales, $\hat\lambda$ and $\hat\tau$, are not considered, and the influences of $\Delta \hx/\hat{\lambda}_0$ and $\Delta \hht/\hat{\tau}_0$ on hydrodynamics cannot be identified.
}

 \vspace{2mm}

\section{Example I: UP properties of DUGKS}
\label{sec:DUGKS}

In this section, we will take the discrete unified gas-kinetic scheme (DUGKS)  \cite{ref:DUGKS15,ref:DUGKS13,ref:DUGKS_Rev} as an example to show how the UP properties are analyzed.

\subsection{Formulation of DUGKS}
The DUGKS is a finite-volume discretization of the kinetic equation \eqref{eq:nonDimensional-BGK} for simulating gas flows in all regimes. For simplicity, we will consider the  one-dimensional case, and the flow is assumed to be isothermal and  smooth without shock discontinuities with a constant relaxation time $\tau$. The computational domain will be divided into a number of uniform cells centered at $x_j$ ($j=1, 2, \cdots, N$) with cell size $\Delta x$, and the interface between cell $j$ and $j+1$ is denote by $x_{j+1/2}$.
The DUGKS can then be expressed as
\begin{equation}
\label{eq:DUGKS}
\dfrac{f_j^{n+1}-f_j^n}{\Delta t} +\xi \dfrac{f_{j+1/2}^{n+1/2}-f_{j-1/2}^{n+1/2}}{\Delta x} =  \dfrac{1}{2\epsilon}\left[Q_j^{n}+Q_j^{n+1}\right],
\end{equation}
where $f_j^n=f_h(x_j, \xi, t_n)$ is the cell-averaged numerical solution over the cell centered at $x_j$ at time $t_n=n\Delta t$, and $Q_j^n = Q_h(x_j, \xi, t_n)=-\tfrac{1}{\tau}[f_h(x_j, \xi, t_n)-f_h^{(eq)}(x_j, \xi, t_n)]$ is the corresponding cell-averaged numerical collision term. Note that here the mid-point and trapezoidal quadrature rules are applied to the convection and collision terms, respectively. The distribution function at cell interface at the half time step $f_{j+1/2}^{n+1/2}=f_h(x_{j+1/2}, \xi, t_n+\Delta t/2)$ is re-constructed by integrating the kinetic equation along the characteristic line with a half time step,
\begin{equation}
\label{eq:DUGKS-flux}
f_{j+1/2}^{n+1/2}-f_{j'}^n =  \dfrac{1}{2\epsilon}\left[Q_{j'}^{n}+Q_{j+1/2}^{n+1/2}\right],
\end{equation}
where $f_{j'}^n = f_h(x_{j+1/2}-\xi\Delta t/2, t_n)$ is the distribution function at the starting point. Note that the implicit discretization in Eqs. \eqref{eq:DUGKS} and \eqref{eq:DUGKS-flux} can be removed by introducing two auxiliary distribution functions in practical computations \cite{ref:DUGKS13}. In DUGKS, $f_{j'}^n$ is linearly interpolated from the cell-center distribution function,
\begin{equation}
f_{j'}^n =  f_{j+1/2}^n -\frac{\Delta t}{2}\xi\sigma_{j+1/2}^n,
\end{equation}
where $\sigma_{j+1/2}^n$ is the slope. For smooth flows, the following approximation can be employed,
\begin{equation}
f_{j+1/2}^n=\dfrac{f_j^n+f_{j+1}^n}{2},\quad \sigma_{j+1/2}^n=\dfrac{f_{j+1}^n-f_j^n}{\Delta x}.
\end{equation}
$Q_{j'}^n$ in \eqref{eq:DUGKS-flux} can be obtained similarly.
Then from Eq. \eqref{eq:DUGKS-flux}, $f_{j+1/2}^{n+1/2}$ can be expressed as
\begin{equation}
\label{eq:f_interface}
f_{j+1/2}^{n+1/2} = \left(\tfrac{1}{2}-\beta\right)\left[f_{j+1}^{n}+\tfrac{\Delta t}{4\epsilon}Q_{j+1}^n\right] + \left(\tfrac{1}{2}+\beta\right)\left[f_{j}^{n}+\tfrac{\Delta t}{4\epsilon}Q_{j}^n\right] + \tfrac{\Delta t}{4\epsilon}Q_{j+1/2}^{n+1/2},
\end{equation}
where $\beta=\tfrac{1}{2}\xi\Delta t  /\Delta x$.

\subsection{Uniform stability}
Regarding the numerical stability of DUGKS, it is not a main task of the present work to give a rigorous analysis. Hence, it is just estimated {\it heuristically} as follows. First, since the collision term in DUGKS (see Eqs. \eqref{eq:DUGKS} and \eqref{eq:f_interface}) is integrated with the trapezoidal rule, which is a semi-implicit time discretization, roughly the small parameter $\epsilon$ in the collision term has no direct influence on the time step.
Indeed, we can think DUGKS as the average of the explicit and the implicit schemes,
\begin{eqnarray}
\frac{f_j^{n+1}-f_j^n}{\Delta t} +\xi \frac{f_{j+\frac 12}^n-f_{j-\frac 12}^n}{\Delta x} =\frac{1}{\epsilon} Q_j^n, \ \ \ \\
\frac{f_j^{n+1}-f_j^n}{\Delta t} +\xi \frac{f_{j+\frac 12}^{n+1}-f_{j-\frac 12}^{n+1}}{\Delta x} =\frac{1}{\epsilon} Q_j^{n+1},\\
f_{j+\frac 12}^{n+\frac 12} = \frac 12(f_{j+\frac 12}^n+f_{j+\frac 12}^{n+1}) +O(\Delta t^2).\ \ \ \ \ \
\end{eqnarray}
Then DUGKS is equivalent to the Crank-Nicolson type discretization of \eqref{eq:BGK} within second order accuracy, being equipped  with the IMEX character.   As is well-known that the Crank-Nicolson type discretizaion is unconditional stable, so is DUGKS as demonstrated by  a number of numerical results in previous studies  \cite{ref:DUGKS15,ref:DUGKS13, ref:DUGKS-A3,ref:DUGKS-A2,ref:DUGKS-A1},  although  the convection CFL constraint,
\begin{equation}
\label{eq:CFL}
\Delta t \le \dfrac{\eta \Delta x}{|\xi|_{\max}},
\end{equation}
is still required theoretically because of the interaction of transport and collision, where $0< \eta < \eta_0$ is the CFL number and $|\xi|_{\max}$ is the maximum discretized particle velocity, $\eta_0$ is some constant.
\vspace{0.2cm}

We just point out that DUGKS has the characteristic of Lax-Wendroff type schemes, thanks to \eqref{eq:DUGKS-flux} or \eqref{eq:f_interface} that has the implication of temporal-spatial coupling.  This is evident in the analysis in Subsection \ref{sub:me} below.

{
\subsection{Collisionless limit}\label{sub:collisionlessLimit}
In the collisionless limit as $\epsilon \to \infty$, from Eq. \eqref{eq:f_interface} we can obtain the cell-interface distribution function at half time step,
\begin{equation}
\label{eq:f-clsInterface}
f_{j+1/2}^{n+1/2} = \dfrac{1}{2}\left(f_{j}^n+f_{j+1}^n\right) - \dfrac{\xi\Delta t}{2\Delta x}\left(f_{j+1}^n-f_{j}^n\right)
\end{equation}
and the DUGKS \eqref{eq:DUGKS} becomes
\begin{equation}
\label{eq:DUGKS-LW}
\dfrac{f_j^{n+1}-f_j^n}{\Delta t} +\dfrac{\xi}{2\Delta x} \left(f_{j+1}^{n}-f_{j-1}^{n} \right) -\dfrac{\xi^2\Delta t}{2\Delta x^2} \left(f_{j+1}^{n}-2f_j^n+f_{j-1}^{n} \right) = 0,
\end{equation}
which is just the Lax-Wendroff scheme for the collisionless BGK equation. This means that the DUGKS is a legitimate scheme in the limit $\epsilon\to \infty$, and the stability condition is
\begin{equation}
\dfrac{|\xi|_{max}\Delta t}{\Delta x}\le 1,
\end{equation}
which is consistent with the CFL condition \eqref{eq:CFL}.
}

\subsection{Modified equation analysis}\label{sub:me}
In order to obtain the modified equation of the DUGKS, we first perform Taylor expansions of $f_{j\pm 1/2}^{n+1/2}$ defined by Eq. \eqref{eq:f_interface} at $(x,t)=(x_j, t_n)$, leading to the following equation after some standard algebraic manipulations,
\begin{eqnarray}
&& \dfrac{f_{j+1/2}^{n+1/2} - f_{j-1/2}^{n+1/2}}{\Delta x} \nonumber \\
&=& \partial_x f_h + \dfrac{\Delta x^2}{6}\partial_x^3 f_h  -\frac{\Delta t}{2}\xi \partial_x^2 f_h  + \frac{\Delta t}{2\epsilon}\partial_x Q_h -\dfrac{1}{2\epsilon}\left(\frac{\Delta t}{2}\right)^2\xi\partial_x^2 Q_h  \nonumber \\
&& + \dfrac{1}{2\epsilon}\left(\frac{\Delta t}{2}\right)^2\partial_x\partial_t Q_h
 + O(\Delta t^3,\Delta t\Delta x^2){\cal L}\left(\tfrac{1}{\epsilon}Q_h\right) +  O(\Delta x^3,\Delta t\Delta x^2){\cal L} f_h,
\end{eqnarray}
where ${\cal L} (\cdot)$ is a general linear operator acting on the dummy variable, representing the collection of temporal and spatial derivatives of different orders. It can take different forms at different places. Here and below we ignore the difference when no confusion is caused. The prefactor in the front of $\cal L$ denotes the corresponding leading orders in $\Delta x$ and $\Delta t$ in the collection.
Note that the collision term also appears in the discretization of the flux, which is a special feature of DUGKS and implies that both the transport and collision effects of particles be precisely included in the scheme. With the above result and performing Taylor expansions of other two terms in Eq. \eqref{eq:DUGKS} at $t_n$, we can obtain that
\begin{eqnarray}
\label{eq:ME-1}
\partial_t f_h + \xi \partial_x f_h  &+&  \dfrac{\Delta t}{2}\underbrace{\left[\partial_t^2 f_h -\xi^2\partial_x^2 f_h -\dfrac{1}{\epsilon}\partial_t Q_h + \dfrac{1}{\epsilon}\xi\partial_x Q_h\right]}_A \nonumber\\
&+& \dfrac{\Delta t^2}{6}\underbrace{\left[\partial_t^3 f_h + \dfrac{3}{4\epsilon}\left(\xi\partial_x\partial_t Q_h -\xi^2\partial_x^2 Q_h-2\partial_t^2 Q_h\right)\right]}_B  + \dfrac{\Delta x^2}{6}\xi\partial_x^3 f_h  \nonumber \\
 &=& \dfrac{1}{\epsilon}Q_h
+ O(\Delta t^3,\Delta t\Delta x^2){\cal L}\left(\tfrac{1}{\epsilon}Q_h\right) +  O(\Delta x^3,\Delta t\Delta x^2, \Delta t^3){\cal L} f_h.
\end{eqnarray}
The high order time derivatives of $f_h$ can be replaced in terms of spatial derivatives successively using the equation \eqref{eq:ME-1}, as done in the modified equation approach \cite{ref:ME}.  We will discuss this below.

For the underbraced term $A$, we can obtain by making use of Eq. \eqref{eq:ME-1} that
\begin{eqnarray}
A&=&(\partial_t-\xi\partial_x)\left(\partial_t f_h+\xi\partial_x f_h-\dfrac{1}{\epsilon}Q_h\right) \nonumber \\
&=& -\dfrac{\Delta t}{2}\left(\partial_t-\xi \partial_x \right)A + O(\Delta t^2){\cal L}\left(\tfrac{1}{\epsilon}Q_h\right) +  O(\Delta x^2,\Delta t^2){\cal L} f_h,
\end{eqnarray}
which gives that
\begin{equation}
A= O(\Delta t^2){\cal L}\left(\tfrac{1}{\epsilon}Q_h\right) +  O(\Delta x^2,\Delta t^2){\cal L} f_h.
\end{equation}
Therefore, the modified equation \eqref{eq:ME-1} can be rewritten as
\begin{equation}
\label{eq:DUGKS-ME2}
\partial_t f_h + \xi \partial_x f_h +\dfrac{\Delta t^2}{6} B + \dfrac{\Delta x^2}{6}\xi\partial_x^3 f_h = \dfrac{1}{\epsilon} Q_h + \mbox{HoT},
\end{equation}
where $\mbox{HoT}$ is the sum of high-order terms,
\begin{equation}
\label{eq:HoT}
\mbox{HoT}=O(\Delta t^3,\Delta t\Delta x^2){\cal L}\left(\tfrac{1}{\epsilon}Q_h\right) +  O(\Delta x^3,\Delta t\Delta x^2, \Delta t^3){\cal L} f_h.
\end{equation}
The first term in HoT contains the scaling parameter $\epsilon$, but we will show below that it has the same order as the second term, and thus HoT depends only on $\Delta t$ and $\Delta x$. To this end, we rewrite Eq. \eqref{eq:DUGKS-ME2} as,
\begin{equation}
\left[1+O(\Delta t^2){\cal L}\right]\left(\tfrac{1}{\epsilon}Q_h\right)=\partial_t f_h + \xi \partial_x f_h +\dfrac{\Delta t^2}{6} \partial_t^3 f_h + \dfrac{\Delta x^2}{6}\xi\partial_x^3 f_h +  O(\Delta x^3,\Delta t\Delta x^2, \Delta t^3){\cal L} f_h,
\end{equation}
suggesting that
\begin{equation}
\label{eq:Q-Df}
\dfrac{Q_h}{\epsilon} = \partial_t f_h + \xi \partial_x f_h + O(\Delta x^2, \Delta t^2) {\cal L} f_h.
\end{equation}
With this estimation, the underbraced term $B$ in Eq. \eqref{eq:ME-1} can be rewritten as
\begin{equation}
B=-\frac{1}{2} \partial_{t}^{3} f_h-\frac{3}{4} \xi^{3} \partial_{x}^{3} f_h-\frac{3}{4} \xi \partial_{x} \partial_{t}^{2} f_h +  O(\Delta x^2, \Delta t^2) {\cal L} f_h,
\end{equation}
and HoT given by \eqref{eq:HoT} is
\begin{equation}
\mbox{HoT}= O(\Delta x^3,\Delta t\Delta x^2, \Delta t^3){\cal L} f_h = O(\Delta x^3,\Delta t\Delta x^2, \Delta t^3),
\end{equation}
which is independent of $\epsilon$. We will neglect ${\cal L} f_h$ in what follows for simplicity, since $\epsilon$ does not appear in these terms. Then the modified equation \eqref{eq:ME-1} can be rewritten as
\begin{equation}
	\label{eq:DUGKS-ME3}
	\partial_t f_h + \xi \partial_x f_h +\dfrac{\Delta t^2}{6} B + \dfrac{\Delta x^2}{6}\xi\partial_x^3 f_h = \dfrac{1}{\epsilon} Q_h + O(\Delta x^3,\Delta t\Delta x^2, \Delta t^3).
\end{equation}


It is clear that as $\Delta t \to 0$ and $\Delta x \to 0$, the modified equation \eqref{eq:DUGKS-ME3} reduces to the kinetic equation \eqref{eq:BGK}, suggesting that the DUGKS consisting of Eqs. \eqref{eq:DUGKS} and \eqref{eq:f_interface} is a consistent second-order scheme in both time and space for the BGK equation \eqref{eq:nonDimensional-BGK}.

\subsection{Analysis of unified properties}

For the asymptotic property of DUGKS at small $\epsilon$, we have the following result.
\begin{theorem}
\label{theorem}
As {{$\Delta t=O(\epsilon^{\alpha})$ and $\Delta x=O(\epsilon^{\beta})$ with $0.5 <\alpha<1$ and $0.5 <\beta<1$}}, the Chapman-Enskog expansion coefficients of $f_h$ obtained from DUGKS satisfy Eq. \eqref{eq:UP} for $n=2$ if the relaxation time $\tau$ is constant.
\end{theorem}

\begin{proof}
With the condition for $\Delta t$ and $\Delta x$, we can write the modified equation \eqref{eq:DUGKS-ME3} as
\begin{equation}
	\label{eq:DUGKS-ME4}
\epsilon\partial_t f_h + \epsilon \xi \partial_x f_h + \dfrac{a}{6}\epsilon^{2\alpha+1} B + \dfrac{b}{6}\epsilon^{2\beta+1}\xi\partial_x^3 f_h =  Q_h + o(\epsilon^2),
\end{equation}
where $a=O(1)$ and $b=O(1)$ are two constants.
Substituting the Chapman-Enskog expansions \eqref{eq:fhCE} and \eqref{eq:tCE} into the above equation leads to the successive equations in terms of the orders of $\epsilon$.
{Note that $2<2\alpha+1<3$ and $2<2\beta+1<3$, suggesting that $\epsilon^{2\alpha+1}=o(\epsilon^2)$, $\epsilon^3=o(\epsilon^{2\alpha+1})$, $\epsilon^{2\beta+1}=o(\epsilon^2)$, and  $\epsilon^3=o(\epsilon^{2\beta+1})$. Then it can be checked straightforwardly that the balance equations at the first three orders are respectively,
\begin{subequations}
\begin{equation}
f_h^{(0)}=f_h^{(eq)},
\end{equation}
\begin{equation}
	D_0f_h^{(0)}=-\dfrac{1}{\tau}f_h^{(1)},
\end{equation}
\begin{equation}
	\partial_{t_1} f_h^{(0)} + D_0f_h^{(1)}=-\dfrac{1}{\tau}f_h^{(2)}.
\end{equation}
\end{subequations}
On the other hand, the terms associated with $\Delta x$ and $\Delta t$ will affect the discrete balance equation at the order of $\epsilon^3$, suggesting  a consistent $f_h^{(3)}$ defined by Eq. \eqref{eq:UP} cannot be achieved. The proof is completed.
}
\end{proof}

In summary, the arguments given in above subsections show that DUGKS is a consistent scheme for the kinetic equation \eqref{eq:BGK} with uniform stability in $\epsilon$, and with Theorem \ref{theorem}, it can be concluded that the DUGKS is a second-order UP scheme, which gives the Navier-Stokes solutions at the cell resolution $\Delta t=o(\sqrt{\epsilon})$ and $\Delta x=o(\sqrt{\epsilon})$.
\vspace{0.2cm}

\begin{remark}
If the collision term is neglected in the reconstruction of the interface distribution function, i.e., $f_{j+1/2}^{n+1/2}=f_{j'}^n = (\tfrac{1}{2}-\beta)f_{j+1}^{n} + (\tfrac{1}{2}+\beta)f_{j}^{n}$, then the modified equation of this collision-less reconstruction (CLR) scheme is (see Appendix A for details)
\begin{eqnarray}
\partial_t f_h &+& \xi \partial_x f_h  + \dfrac{\Delta x^2}{6}\xi\partial_x^3 f_h + \dfrac{\Delta t}{2}\underbrace{\left[\partial_t^2 f_h -\xi^2\partial_x^2 f_h -\dfrac{1}{\epsilon}
\partial_t Q_h\right]}_{A'}\nonumber \\
 &=&\dfrac{1}{\epsilon}Q_h + O(\Delta t^2){\cal L} \left(\tfrac{1}{\epsilon}Q_h\right)+ O(\Delta x^3, \Delta t\Delta x^2, \Delta t^2){\cal L} f_h .
\end{eqnarray}
It can be shown that $A' = -\xi\partial_x\partial_t f_h - \xi^2 \partial_x^2 f_h + O(\Delta t, \Delta x^2){\cal L} f_h$, and the modified equation can be rewritten as
\begin{equation}
\partial_t f_h + \xi \partial_x f_h  + \dfrac{\Delta x^2}{6}\xi\partial_x^3 f_h - \dfrac{\Delta t}{2}\left(\xi\partial_x\partial_t f_h + \xi^2 \partial_x^2 f_h\right) = \dfrac{1}{\epsilon}Q_h  + O(\Delta x^3, \Delta t^2),
\end{equation}
which suggests that the UP order of the scheme would degenerate in comparison with the DUGKS due to the non-vanishing $A'$ (see Eq. \eqref{eq:DUGKS-ME2}). As a result, under the same cell resolution the above scheme cannot give accurate Navier-Stokes solutions. This result confirms that it is important to consider the collision term in the reconstruction of numerical flux \cite{ref:ChenXu-2015}, and implies the essence of the temporal-spatial coupling in the design of schemes \cite{ref:Li-TS-2019,ref:Li-Du-2016, ref:Pan-Xu-Li-Li-2016}.
\end{remark}

\subsection{Numerical test}
\label{subsec-num-1}
We now test the UP property of the DUGKS with the two-dimensional incompressible Taylor vortex in a periodic domain $0\le x \le 1$ and $0\le y \le 1$. At the hydrodynamic scale ($\epsilon\ll 1$), the flow is governed by the incompressible Navier-Stokes equations and has the following analytical solution,
\begin{subequations}
\begin{equation}
u_x=-\dfrac{u_0}{A}\cos(Ax)\sin(By)e^{-\nu\theta t},
\end{equation}
\begin{equation}
u_y=\dfrac{u_0}{B}\sin(Ax)\cos(By)e^{-\nu\theta t},
\end{equation}
\begin{equation}
p(x, y, t)=p_0-\dfrac{\rho_0 u_{0}^{2}}{4}\left[\frac{\cos (2 A x)}{A^{2}}+\frac{\cos (2 B y)}{B^{2}}\right] e^{-2 \nu \alpha t},
\end{equation}
\end{subequations}
where $u_0$ is a constant, $\theta = A^2 + B^2$, $\nu$ is the shear viscosity, and $\u = (u_x,u_y)$ and $p$ are the velocity and pressure, respectively, $p_0=\rho_0 RT_0$ is the reference pressure with $\rho_0$ the average density and $T_0$ the constant temperature.

For this low-speed isothermal flow, the discrete velocity set used in the DUGKS is chosen based on the three-point Gauss-Hermite quadrature in each direction as shown in \cite{ref:DUGKS13}, namely, $\bxi_0=(0,0)$, $\bxi_1=-\bxi_3=c(1, 0)$, $\bxi_2=-\bxi_4=c(0, 1)$, $\bxi_5=-\bxi_7=c(1, 1)$, $\bxi_6=-\bxi_8=c(-1, 1)$, with $c=\sqrt{3RT_0}$. The equilibrium distribution function is approximated with the low Mach number expansion of the Maxwellian equilibrium,
\begin{equation}
f_i^{(eq)}=w_i\rho\left[1+\dfrac{\bxi_i \cdot \u}{R T_0}+\dfrac{(\bxi_i \cdot \u)^{2}}{2(R T_0)^{2}}-\dfrac{|\u|^{2}}{2 R T_0}\right],
\end{equation}
where $w_0=4/9$, $w_1=w_2=w_3=w_4=1/9$, and $w_5=w_6=w_7=w_8=1/36$. It is clear that these parameters are the same as the standard D2Q9 lattice Boltzmann equation model \cite{ref:GuoBook}.

In our simulations, we set $A = B = 2\pi$, $u_0 = 0.01$, and $RT_0 = 0.5$, such that the Mach number is small and the flow can be well recognized as incompressible. The relaxation time is determined from the shear viscosity, $\tau=\nu/RT_0$ so that the parameter $\epsilon$ can be adjusted by changing the value of $\nu$. Uniform meshes are employed and the CFL number $\eta$ is set to be 0.5 for each mesh. Periodic boundary conditions are imposed on all boundaries, and the distribution functions are initialized by setting $f_i=f_i^{(eq)}+\epsilon f_i^{(1)}$, which is the Chapman-Enskog approximation at the Navier-Stokes order.

\begin{figure}[!htb]
\centering
\includegraphics[width=0.6\textwidth]{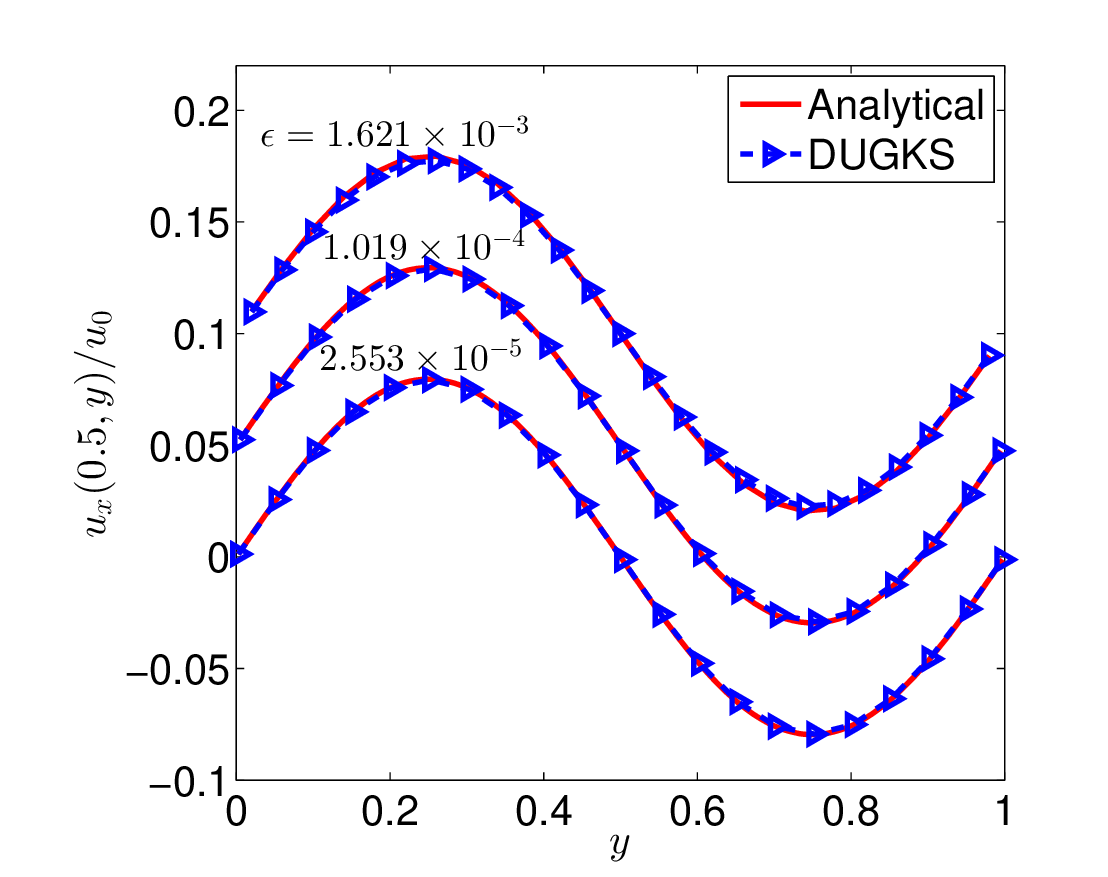}
\caption{Velocity profiles of the Taylor vortex flow at $t=t_c$ predicted by the DUGKS with $\Delta x = \sqrt{\epsilon}$. The profiles for $\epsilon=1.621\times 10^{-3}$ and $1.019\times 10^{-4}$ are shifted upward for clarity.}
\label{fig:Taylor}
\end{figure}

Three values of $\epsilon$ for continuum flow regime, i.e., $1.621\times 10^{-3}$, $1.019\times 10^{-4}$, and $2.553\times 10^{-5}$, are considered in the simulations. We first test whether the Navier-Stokes solution for each case can be captured by the DUGKS with a uniform mesh with resolution of $\Delta x \sim \epsilon^{\alpha}$ with $\alpha=0.501$, (and thus $\Delta t\sim\epsilon^{\alpha}$). Specifically, uniform meshes with size of $25\times 25$, $100\times 100$, and $200\times 200$ are adopted for $\epsilon=1.621\times 10^{-3}$, $1.019\times 10^{-4}$, and $2.553\times 10^{-5}$, respectively. The velocity profiles at $t=t_c\equiv \ln2/(\nu\theta)$, at which the magnitude of the velocity decays to one half of the original one, are measured and shown in Fig. \ref{fig:Taylor}. It can be observed that the velocity profiles are well captured by the DUGKS with the meshes, confirming its second-order UP property. We note that the mesh resolutions are much larger than the kinetic scale in the tests ($\Delta x/\epsilon=24.68$, 98.18, and 195.81 respectively for the three cases).

\begin{figure}[!htb]
\centering
\includegraphics[width=0.6\textwidth]{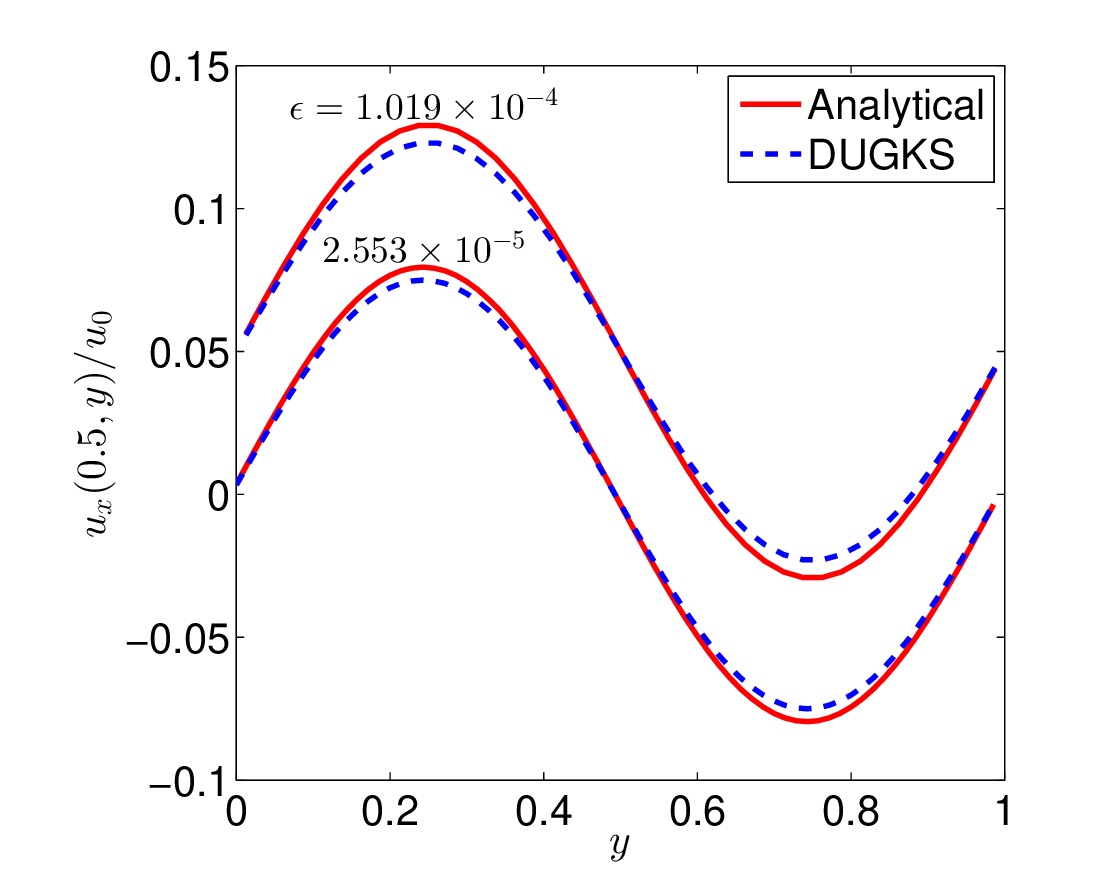}
\caption{Velocity profiles of the Taylor vortex flow at $t=t_c$ predicted by the DUGKS with $\Delta x \approx \epsilon^{0.4}$. The profile for $1.019\times 10^{-4}$ is shifted upward for clarity.}
\label{fig:Taylor_coarse}
\end{figure}

We then tests whether the DUGKS can capture the Navier-Stokes solutions with meshes of coarser resolution than $O(\epsilon^{0.5})$.  As an example, we use a $40\times 40$ mesh and a $70\times 70 $ one for  $\epsilon=1.019\times 10^{-4}$ and $2.553\times 10^{-5}$, respectively, such that $\Delta x \approx \epsilon^{0.4}$. The velocity profiles are shown in Fig. \ref{fig:Taylor_coarse}, and some clear deviations between the numerical and analytical solutions can be observed.
These results confirm the UP analysis of DUGKS presented above, namely, $\Delta x$ and $\Delta t$ should be of $o(\epsilon^{0.5})$ in order to capture the Navier-Stokes hydrodynamics.

\begin{figure}[!htb]
\centering
\includegraphics[width=0.6\textwidth]{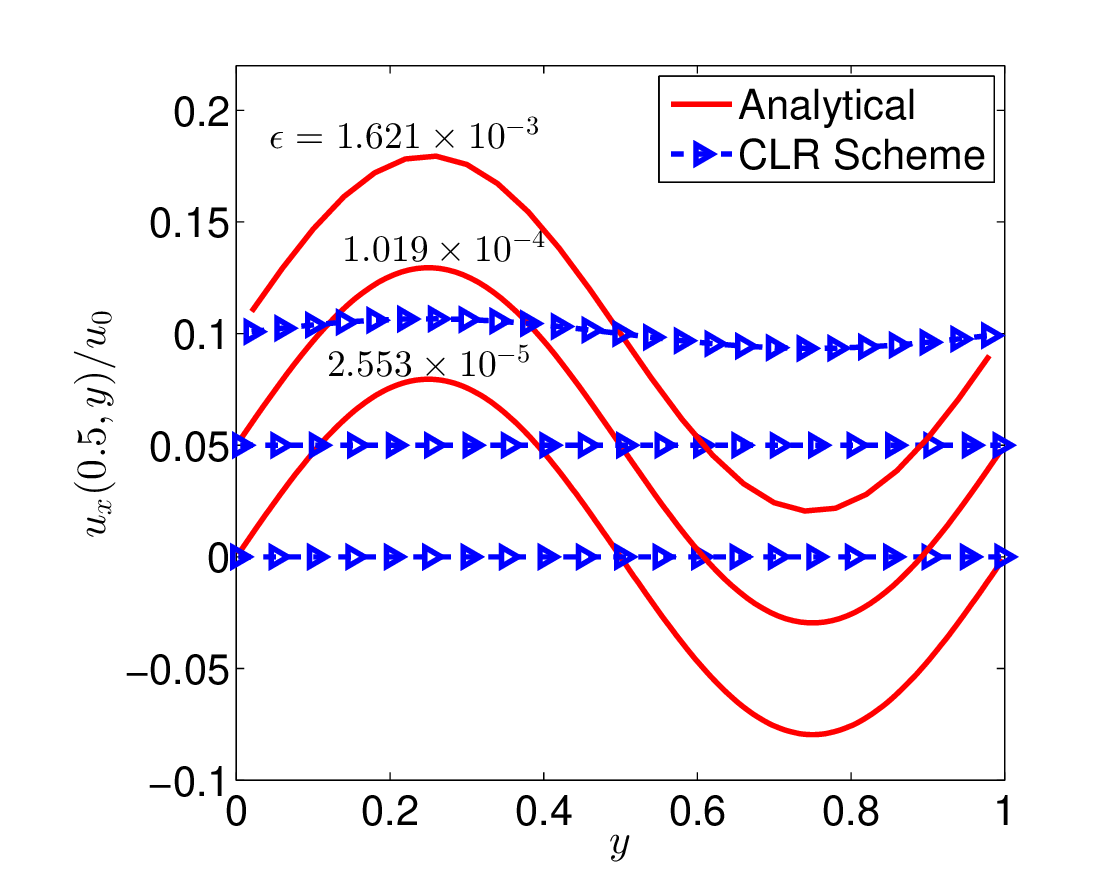}
\caption{Velocity profiles of the Taylor vortex flow at $t=t_c$ predicted by the CLR scheme with $\Delta x = \epsilon^{0.501}$. The profiles for $\epsilon=1.621\times 10^{-3}$ and $1.019\times 10^{-4}$ are shifted upward for clarity.}
\label{fig:CLR}
\end{figure}

We now test whether the CLR scheme, where the distribution function $f_{j+1/2}^{n+1/2}$ is reconstructed by solving the collision-less kinetic equation as noted in Remark II and Appendix A, can capture the Navier-Stokes solutions with the same numerical resolution as used in the DUGKS (i.e., $\Delta x \sim \epsilon^{0.501}$ with CFL number $\eta=0.5$). The predicated velocity profiles for $\epsilon=1.621\times 10^{-3}$, $1.019\times 10^{-4}$, and $2.553\times 10^{-5}$ are shown in Fig. \ref{fig:CLR}, which clearly demonstrates that the CLR scheme is too dissipative to capture the Navier-Stokes solutions under this cell resolution. These results confirm the UP property of this scheme shown in Appendix.

{
\section{Example II: UP properties of a second-order IMEX Runge-Kutta scheme}
\label{sec:IMEX-RK}
In this section we investigate the UP property of a second-order IMEX Runge-Kutta (IMEX-RK) scheme presented in \cite{ref:AP_NS_Hu2017}. Originally, the time asymptotics of the IMEX-RK scheme was analyzed in  \cite{ref:AP_NS_Hu2017}, in which the spatial gradient was considered to be exact. It was claimed that the scheme could capture the Navier-Stokes limit without resolving the kinetic scheme, i.e., as $\epsilon = o(\Delta t)$ and $\Delta t = o(\epsilon^{0.5})$. However, as we demonstrated in the introduction and Sec. \ref{sec:DUGKS}, it is critical to consider the temporal-spatial coupling for the asymptotic behaviors of a multiscale kinetic scheme. We will revisit this scheme in the UP framework to study how the spatial discretization influences the asymptotics.

\subsection{Formulation}
The second-order IMEX-RK scheme considered here is a two-stage globally stiffly accurate (GSA) one, which is named as IMEX-II-GSA(2,3,2) in Ref. \cite{ref:AP_NS_Hu2017}. In one dimension, the scheme can be expressed as
\begin{subequations}
\label{eq:IMEX-RK}
\begin{equation}
f_h^{n+1/2}=f_h^{n}-\frac{\xi\Delta t}{2} \partial_{x} f_h^{n}+\frac{\Delta t}{2 \epsilon} Q_h^{n+1/2},
\end{equation}
\begin{equation}
f_h^{n+1}=f_h^{n}-\Delta t\xi \partial_{x} f_h^{n+1/2}+\frac{\Delta t}{2 \epsilon}\left[Q_h^{n}+Q_h^{n+1}\right],
\end{equation}
\end{subequations}
where the spatial derivatives are not discretized. Note that the hydrodynamic variables appearing in $f^{(eq)}$ of $Q^{n+1/2}$ and $Q^{n+1}$ are first obtained from \eqref{eq:IMEX-RK} by taking the conservative moments,
\begin{subequations}
	\label{eq:IMEX-W}
	\begin{equation}
		{\bf \bar{W}}_h^{n+1/2}={\bf \bar{W}}_h^n - \frac{\Delta t}{2} \partial_{x} \bm{\mathcal F}_h^{n},
	\end{equation}
	\begin{equation}
		{\bf \bar{W}}_h^{n+1}={\bf \bar{W}}_h^n - \frac{\Delta t}{2} \partial_{x} \bm{\mathcal F}_h^{n+1/2},
	\end{equation}
\end{subequations}
where ${\bf \bar{W}}_h=(\rho_h,\rho_h\u_h, \rho_h E_h)^T$ and $\bm{\mathcal F}_h=\int{\bxi \bm{\psi} f_h d\bxi}$.

\subsection{Time asymptotics}
From Eq. \eqref{eq:IMEX-RK}, we can obtain that
\begin{equation}
\dfrac{f_h^{n+1}-f_h^{n}}{\Delta t}+\xi \partial_{x} f_h^n -\dfrac{\Delta t}{2}\xi^2\partial_x^2 f_h^n+\dfrac{\Delta t}{2\epsilon}\partial_x Q_h^{n+1/2}=\frac{1}{2 \epsilon}\left[Q_h^{n}+Q_h^{n+1}\right].
\end{equation}
By performing the Taylor expansions at $t_n$, we can obtain the modified equation of the IMEX-RK sheme,
\begin{eqnarray}
	\label{eq:ME-RK-1}
	\partial_t f_h + \xi \partial_x f_h  &+&  \dfrac{\Delta t}{2}\underbrace{\left[\partial_t^2 f_h -\xi^2\partial_x^2 f_h -\dfrac{1}{\epsilon}\partial_t Q_h + \dfrac{1}{\epsilon}\xi\partial_x Q_h\right]}_A \nonumber\\
	&+& \dfrac{\Delta t^2}{6}\underbrace{\left[\partial_t^3 f_h + \dfrac{3}{2\epsilon}\left(\xi\partial_x\partial_t Q_h -\partial_t^2 Q_h\right)\right]}_B  + \dfrac{\Delta x^2}{6}\xi\partial_x^3 f_h  \nonumber \\
	&=& \dfrac{1}{\epsilon}Q_h
	+ O(\Delta t^3,\Delta t\Delta x^2){\cal L}\left(\tfrac{1}{\epsilon}Q_h\right) +  O(\Delta x^3,\Delta t\Delta x^2, \Delta t^3){\cal L} f_h,
\end{eqnarray}
which is similar to the modified equation \eqref{eq:ME-1} for the DUGKS. Actually, the underbraced term $A$ is the same as that in Eq. \eqref{eq:ME-1}. Then, after some similar analysis we can obtain the simplified modeified equation for the IMEX-RK scheme,
\begin{equation}
\label{eq:ME-RK-2}
\partial_t f_h + \xi \partial_x f_h + \dfrac{\Delta t^2}{6} B = \dfrac{1}{\epsilon}Q_h + O(\Delta x^3,\Delta t\Delta x^2, \Delta t^3){\cal L} f_h,
\end{equation}
with
\begin{equation}
B=-\dfrac{1}{2}\partial_t^3 f_h + \dfrac{3}{2}\xi^2\partial_x^2\partial_t f_h.
\end{equation}
It is clear that the IMEX-RK is a second-order time discretization of the kinetic equation. Furthermore, it can be shown that if $\Delta t=O(\epsilon^{\alpha})$ with $0.5<\alpha<1$, the IMEX-RK is second-order UP in time. This suggests that it can preserve the Navier-Stokes solution if $\Delta t=o(\epsilon^{0.5})$ and $\epsilon=o(\Delta t)$, which is consistent with the result in \cite{ref:AP_NS_Hu2017}.

It should be noted, however, that the above analysis only considers the temporal discretization, and the spatial gradient is assumped to be exact. Since the kinetic equation is temporal-spatial coupled intrinsickly, the temporal-spatial coupling at discrete level is also critical for a good practical numerical scheme \cite{ref:Li-TS-2019,ref:Li-Du-2016,ref:Pan-Xu-Li-Li-2016}. This is also evident from the latter stability analysis and numerical results.

\subsection{Uniform stability}
The numerical stability of the IMEX-RK depends on the spatial discretization. If the space is discretized with a uniform mesh and the convection term is discretized with the second-order central difference scheme, the IMEX-RK scheme \eqref{eq:IMEX-RK} can be expresed as
\begin{subequations}
	\label{eq:IMEX-RK-space}
	\begin{equation}
		f_j^{n+1/2}=f_j^{n}-\frac{\xi\Delta t}{2\Delta x} \left(f_{j+1/2}^n - f_{j+1/2}^n\right) + \frac{\Delta t}{2 \epsilon} Q_j^{n+1/2},
	\end{equation}
	\begin{equation}
		f_j^{n+1}=f_j^{n}-\dfrac{\Delta t}{\Delta x}\left(f_{j+1/2}^{n+1/2} - f_{j+1/2}^{n+1/2}\right) +\frac{\Delta t}{2 \epsilon}\left[Q_j^{n}+Q_j^{n+1}\right],
	\end{equation}
\end{subequations}
where $f_j^n=f_h(x_j,t_n)$ and
\begin{equation}
 f_{j+1/2}^n=\dfrac{1}{2}\left(f_j^n+f_{j+1}^n\right),\quad f_{j+1/2}^{n+1/2} = \dfrac{1}{2}\left(f_j^{n+1/2} + f_{j+1}^{n+1/2}\right),
\end{equation}
which is similar with the linear reconstruction used in the DUGKS.

In the collisionless limit ($\epsilon\to \infty$), the IMEX-RK with the above spatial discretization becomes
\begin{equation}
	\label{eq:IMEX-RK-LW}
	\dfrac{f_j^{n+1}-f_j^{n}}{\Delta t}+\dfrac{\xi}{2\Delta x}\left(f_{j+1}^n-f_{j-1}^n\right) -\dfrac{\xi^2\Delta t}{8\Delta x^2}\left(f_{j+2}^n-2f_j^n+f_{j-2}^n \right) = 0,
\end{equation}
which is similar with the case of the DUGKS, Eq. \eqref{eq:DUGKS-LW}, except for the second term on the left hand side. Generally this term  introduces artifical dissipation that can enhance numerical stability. However, it can be shown via the von Neumann analysis that the scheme \eqref{eq:IMEX-RK-LW} is {\it unconditionally unstable}. With a finite value of $\epsilon$, the numerical stability of the IMEX-RK can be enhanced, but it is a nontrivial task to prove its stability property rigorously.

The above analysis indicates that the IMEX-RK with the second order central difference discretization is not a legitimate scheme of the kinetic equation in the collisionless limit, and is not uniformly stable in terms of $\epsilon$. Therefore, the IMEX-RK with such discretization is not a UP scheme. Particularly, this scheme cannot be applied to multiscale flows with a wide range of $\epsilon$. Other types of spatial discretization methods can be employed to enhance the numerical stability. For example, a second-order finite-volume scheme with a slope limiter reconstruction was used in \cite{ref:AP-NS-FJJCP2010}, and a fifth order weighted essentially non-oscillatory (WENO) finite-difference method was employed in \cite{ref:AP_NS_Hu2017}. However, these methods may introduce additional numerical dissipation that destroys the UP property.

\subsection{Numerical test}
To test the capability of the IMEX-RK scheme in capturing the Navier-Stokes behavior, we now simulate the Taylor vortex problem with the same parameters as used in Subsec. \ref{subsec-num-1}. In order to overcome the numerical instability problem of the central difference discretization, an improved third-order WENO interpolation \cite{ref:ESWENO}, with better numerical stability and accuracy than the classical WENO method, is employed to calculate $f_{j+1/2}$ in Eq. \eqref{eq:IMEX-RK-space}.

The velocity profiles at $t_c$ for the three Knudsen numbers are shown in Fig. \ref{fig:IMEX-RK}, with $\Delta x\approx \epsilon^{0.501}$. It can be seen that the numerical results predicted by the IMEX-RK scheme agree well with the analytical solutions for $\epsilon=1.019\times 10^{-4}$ and $2.553\times 10^{-5}$, but deviation is clearly observed for the case of $\epsilon=1.621\times 10^{-3}$, which can be attributed to the relatively large numerical dissipation of the WENO interapolation near extremum points.

\begin{figure}[!htb]
	\centering
	\includegraphics[width=0.6\textwidth]{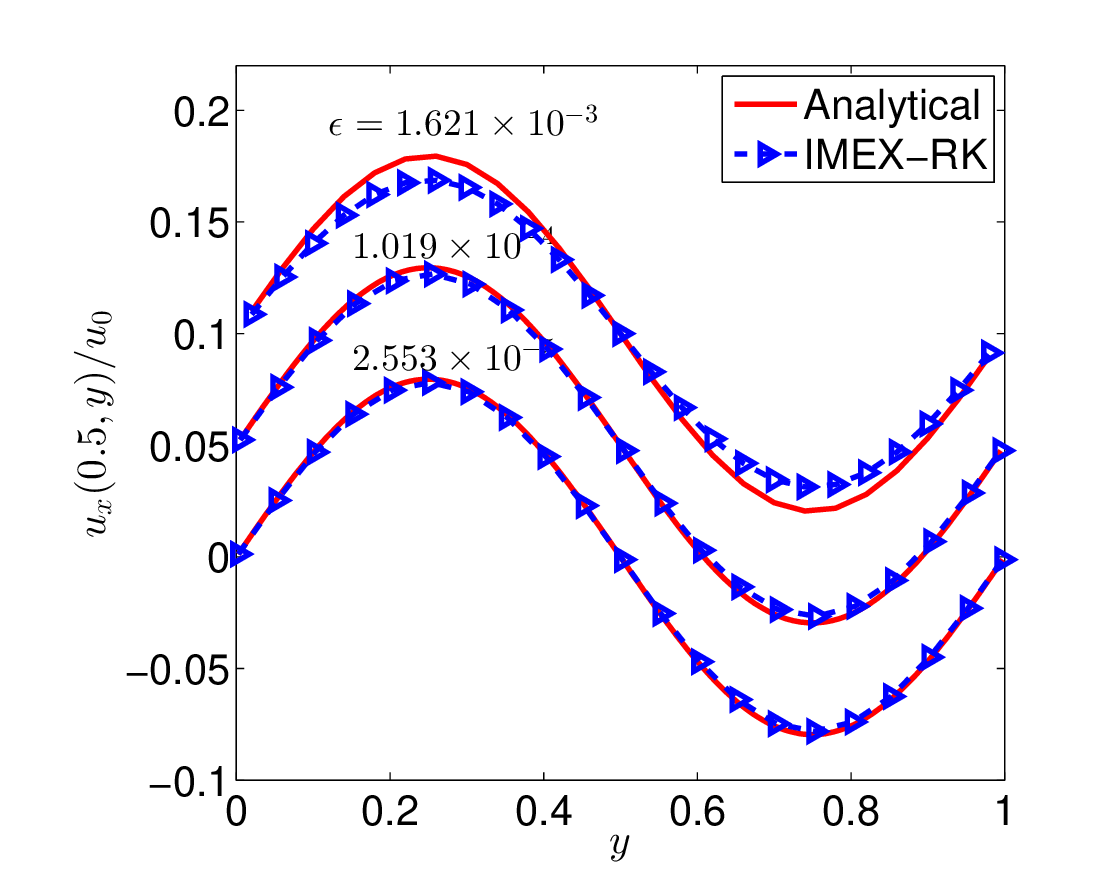}
	\caption{Velocity profiles of the Taylor vortex flow at $t=t_c$ predicted by the IMEX-RK scheme with the improved third-order WENO interpolation and $\Delta x = \epsilon^{0.501}$. The profiles for $\epsilon=1.621\times 10^{-3}$ and $1.019\times 10^{-4}$ are shifted upward for clarity.}
	\label{fig:IMEX-RK}
\end{figure}

}

\section{Summary and Discussions}
\label{sec:summary}

For a kinetic scheme for multi-scale flow simulations, it is required to capture the hydrodynamic behaviors in the continuum limit without resolving the kinetic scales, but with the cell resolution for distinguishing the hydrodynamic structure only. The asymptotic preserving (AP) concept is helpful for understanding the limiting behavior of a kinetic scheme as $\epsilon \to 0$, i.e., the dynamics described by the Euler equations. In the present work, we propose the unified preserving (UP) concept to distinguish the asymptotic behaviors of a kinetic scheme at small but nonzero $\epsilon$ and the order of asymptotics to capture the corresponding hydrodynamic structure.

{
 We note that the role of the UP analysis is different from the general numerical analysis method for kinetic schemes aiming to solve certain hydrodynamic equations, such as the lattice Boltzmann equation and the gas-kinetic scheme. For such numerical methods, the target macroscopic (Euler or Navier-Stokes) equations to be solved are clearly specified, and the classical numerical analysis can provide useful information on how to choose the numerical parameters; on the other hand, the UP analysis is used mainly for kinetic schemes for simulating multiscale flows rather than for certain specific hydrodynamic equations. Actually, such a scheme can be viewed as a dynamic system on a finite mesh with discrete time, and its dynamic behavior depends on the numerical parameters (mesh size and time step) in addition to physical ones. Owing to the multiscale nature of the system, it is not clear what continuum hydrodynamic equations are solved by the kinetic scheme, or even no well-posed hydrodynamics models are available. The goal of the UP analysis is just to identify the hydrodynamic behaviors under different numerical parameters.
}

Generally, a UP scheme with order $n\ge 1$ is also AP, and therefore the UP concept can be viewed as an extension of AP concept. But it should be noted that the approach for analyzing the UP properties is different from that for analyzing the AP property of a kinetic scheme. In the AP analysis, the corresponding discrete formulation of the kinetic scheme as $\epsilon\to 0$ is first derived, and then is proved to be a consistent and stable discretization of the corresponding hydrodynamic equations. On the other hand, in the UP framework the analysis is based on the Chapman-Enskog expansion of the modified equation of the kinetic scheme, and the expansion coefficients at different orders of $\epsilon$ can be compared with those of the original Chapman-Enskog expansion of the kinetic equation and thus the asymptotic degree of accuracy can be assessed.

\begin{figure}[!htb]
	\centering
	\includegraphics[width=0.6\textwidth]{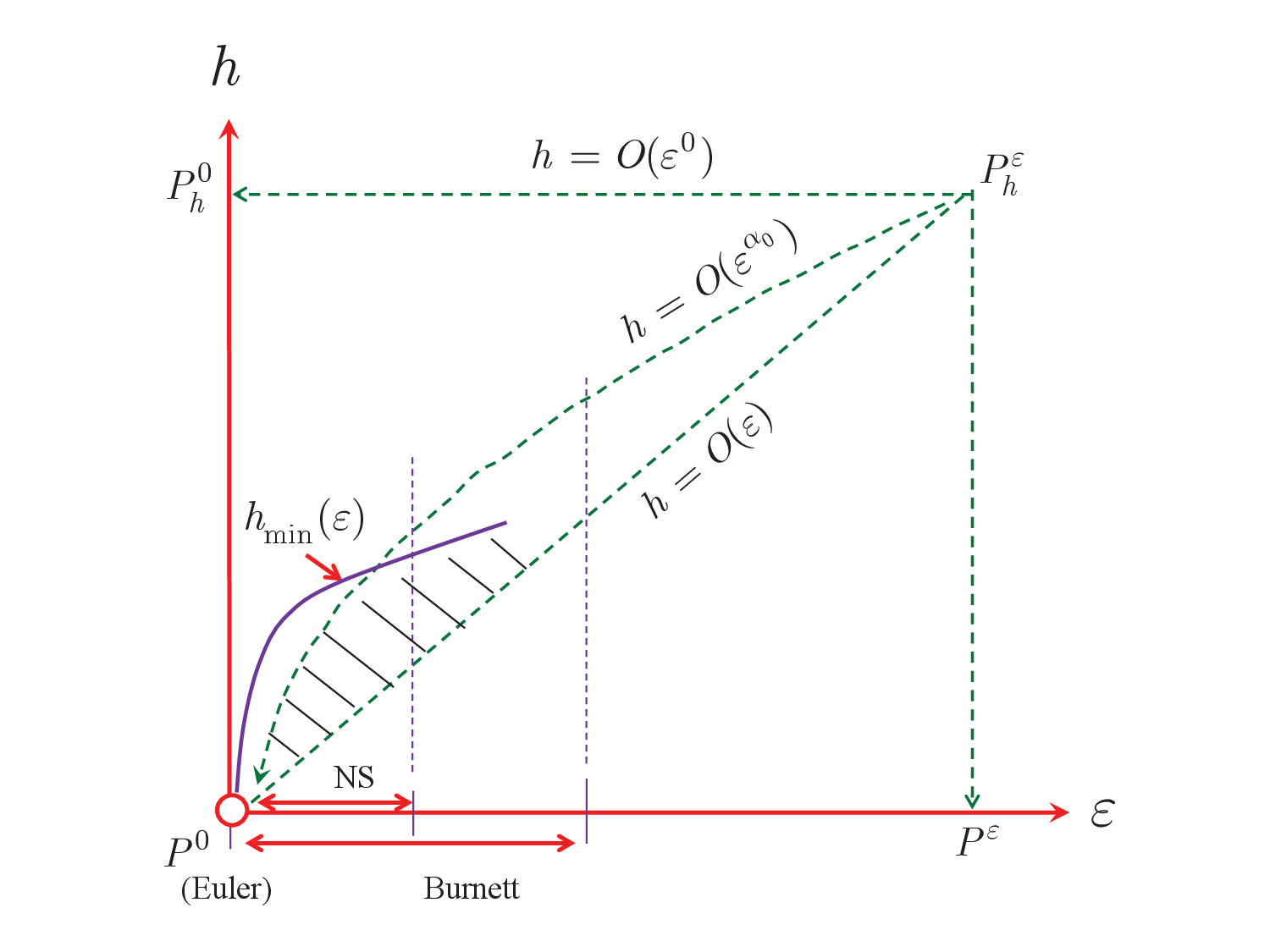}
	\caption{Schematic of asymptotic path to the limiting hydrodynamic regimes with flow structures. The lines are similar with those shown in Fig. \ref{fig:upPath}, and the line $h_{\min}(\epsilon)$ represents the minimum resolution required by flow structures. The shadowed region represents the parameter space for a UP scheme.}
	\label{fig:upPath_h}
\end{figure}

The UP order of a kinetic scheme depends on three scales, i.e., kinetic scale $(\hat{\lambda}_0, \hat{\tau}_0)$, numerical scale $(\Delta \hat{x}, \Delta \hht)$, and hydrodynamic scale $(\hat{l}_0, \hht_0)$. Specifically, the order is related to the dimensionless parameters $\epsilon=\htau_0/\hht_0=\hat{\lambda}_0/\hat{l}_0$, $\delta_x=\Delta\hat{x}/\hat{\lambda}_0$, and $\delta_t=\Delta\hht/\htau_0$. Furthermore, the accuracies of spatial and temporal discretizations also affect the UP order. In general, for given cell size and time step scalings, the UP order increases with the accuracy order of the scheme. It is emphasized here that the present work assumes the numerical resolution is adequate to resolve the flow physics at the hydrodynamic limit. For practical problems with flow structures (such as boundary layer, shear layler, shock layer, and oscillation period) some restrictions to the mesh size and/or time step are required. For instance, in order to capture the shock structure, whose thickness is of $O(\epsilon)$ \cite{ref:BirdBook}, a mesh size of $O(\epsilon)$ is necessary. Under such circumstance, the numerical resolution $h$ should be not larger than an $\epsilon$-dependent minimum $h_{\min}(\epsilon)$, which is an additional imposition on the UP scheme. Consequently, the parameter space of the kinetic scheme to approach its asymptotic limit may be changed, as demonstrated in Fig. \ref{fig:upPath_h}. Only as the numerical resolution of the scheme falls in the shadowed region, the kinetic scheme can capture the correct hydrodynamic behaviors without resolving the kinetic scale. The main target for the UP analysis is to evaluate the kinetic schemes if they can recover the shock capturing schemes for the NS solutions in the hydrodynamic flow regime, such as the laminar boundary layer,  without imposing kinetic scale resolution on the mesh size and time step. 

As an example, the UP property of a temporal-spatial discretized kinetic scheme, DUGKS, is first analyzed. It is shown that the DUGKS is a second-order UP scheme which can capture hydrodynamics at the Navier-Stokes level without resolving the kinetic scale, and the numerical test confirms the property. On the other hand, if the distribution function at cell interface is reconstructed by solving the collision-less kinetic equation, the scheme is only of first-order UP under the same numerical resolution. The results confirm the necessity of the inclusion of collision effect in the reconstruction. We note that some kinetic schemes also consider the collision effect in flux reconstruction \cite{ref:JiangDW,ref:UGKS}, and it would be interesting to analyze the UP properties of these schemes.
{
The UP property of a second-order time-discrete scheme, IMEX-RK, is further analyzed. It is shown the scheme is of second-order time UP provided the spatial gradient is exact. However, if the second-order central-difference (i.e., linear reconstruction) is used as in the DUGKS, the IMEX-RK is not uniformly stable, and thus is not UP. The use of a third-order improved WENO interpolation can enhance the numerical stability, but numerical results demonstrate that the numerical dissipation arising from spatial discretization can destroy the UP property. These findings again confirm the necessarity of considering the fully temporal-spatial particle transport and collision coupling in designing a UP kinetic scheme.
}

The UP property discussed in this study is focused on the case of hydrodynamic regime where $\epsilon$ is small, which is the key requirement for kinetic schemes. For flows beyond the continuum regime, the solution of a UP scheme is expected to be able to capture the non-continuum flow physics since the scheme is consistent with the kinetic equation. For instance, the DUGKS has been successfully applied to a variety of non-equilibrium flows ranging from slip to free-molecular regimes \cite{ref:GDVM-WP,ref:DUGKS-A3,ref:DUGKS-A2,ref:DUGKS-A1}.

Finally, we remark that other techniques rather than the original Chapman-Enskog expansion can also be employed to analyze the asymptotic property of a kinetic scheme in the UP framework, such as the {the regularized Chapman-Enskog expansions \cite{ref:RegCE_Rose1989, ref:NormCE_Slemrod1989,ref:SumCE_Gorban1982}, the Maxwell iteration procedure \cite{ref:HarrisBook}, or the order-of-magnitude method \cite{ref:Struchtrup2004,ref:Kauf2010scale}}.

\section*{Acknowledgments}
ZLG acknowledges the support by the National Natural Science Foundation of China (51836003 and 11872024). JQL was supported by Natural Science Foundation of China (11771054, 12072042 and 91852207) and the Sino-German Research Group Project (No. GZ1465). KX was supported by Hong Kong Research Grant Council (16208021) and National Natural Science Foundation of China (11772281 and 91852114). Helpful discussions with Prof. Songze Chen, Prof. Lei Wu, and Dr. Chang Liu are gratefully acknowledged.

\appendix
\section{UP analysis of the CLR scheme}
\label{sec:appen}
In this appendix we will derive the modified equation of the CLR kinetic scheme noted in Remark II in which the cell-interface distribution function $f_{j+1/2}^{n+1/2}$ is reconstructed from the collision-less BGK equation,
\begin{equation}
\label{eq:sf_interface}
f_{j+1/2}^{n+1/2} = f_{j'}^n = \left(\tfrac{1}{2}-\beta\right)f_{j+1}^{n} + \left(\tfrac{1}{2}+\beta\right)f_{j}^{n}.
\end{equation}
Then we have
\begin{equation}
\dfrac{f_{j+1/2}^{n+1/2} - f_{j-1/2}^{n+1/2}}{\Delta x}
 =\partial_x f_h + \dfrac{\Delta x^2}{6}\partial_x^3 f_h  -\frac{\Delta t}{2}\xi \partial_x^2 f_h  + O(\Delta x^3, \Delta t\Delta x^2){\cal L} f_h.
\end{equation}
The modified equation can thus be derived as
\begin{equation}
\label{eq:sME-1}
\begin{array}{rl}
&\partial_t f_h + \xi \partial_x f_h  + \dfrac{\Delta t}{2}\underbrace{\left[\partial_t^2 f_h -\xi^2\partial_x^2 f_h -\dfrac{1}{\epsilon}
\partial_t Q_h\right]}_{A'}  + \dfrac{\Delta x^2}{6}\xi\partial_x^3 f_h  \\
=& \dfrac{1}{\epsilon}Q_h
+ O(\Delta t^2){\cal L} \left(\tfrac{1}{\epsilon}Q_h\right)+ O(\Delta x^3, \Delta t\Delta x^2, \Delta t^2){\cal L} f_h.
\end{array}
\end{equation}
With this equation we can estimate the underbraced term $A'$,
\begin{eqnarray}
A'&=&\left(\partial_t-\xi\partial_x\right)\left(\partial_t f_h +\xi\partial_x f_h -\dfrac{1}{\epsilon}Q_h\right) - \dfrac{1}{\epsilon}\xi\partial_x Q_h\nonumber \\
&=& \dfrac{\Delta t}{2}\left(\xi \partial_x  - \partial_t  \right)A'  -\dfrac{1}{\epsilon}\xi\partial_x Q_h
+ O(\Delta t^2){\cal L} \left(\tfrac{1}{\epsilon}Q_h\right)+ O(\Delta x^2, \Delta t^2){\cal L} f_h,
\end{eqnarray}
which leads to
\begin{equation}
A'= -\dfrac{1}{\epsilon}\xi\partial_x Q_h + O(\Delta t){\cal L} \left(\tfrac{1}{\epsilon}Q_h\right)+ O(\Delta x^2, \Delta t^2){\cal L} f_h.
\end{equation}
Therefore, the modified equation \eqref{eq:sME-1} can be rewritten as
\begin{eqnarray}
\label{eq:sME-2}
\partial_t f_h &+& \xi \partial_x f_h   - \dfrac{\Delta t}{2\epsilon}\xi\partial_x Q_h + \dfrac{\Delta x^2}{6}\xi\partial_x^3 f_h \nonumber \\
&=& \dfrac{1}{\epsilon}Q_h  + O(\Delta t^2){\cal L} \left(\tfrac{1}{\epsilon}Q_h\right)+ O(\Delta x^3, \Delta t\Delta x^2, \Delta t^3){\cal L} f_h,
\end{eqnarray}
from which we can obtain that
\begin{equation}
\dfrac{1}{\epsilon}Q_h = \partial_t f_h + \xi \partial_x f_h + O(\Delta t, \Delta x^2){\cal L} f_h.
\end{equation}
Thus,
\begin{equation}
A'= -\dfrac{1}{\epsilon}\xi\partial_x Q_h + O( \Delta t, \Delta x^2){\cal L} f_h = -\xi\partial_x\partial_t f_h - \xi^2 \partial_x^2 f_h + O(\Delta t, \Delta x^2){\cal L} f_h,
\end{equation}
and the modified equation \eqref{eq:sME-2} can be re-expressed as
\begin{equation}
\label{eq:sME-2x}
\partial_t f_h + \xi \partial_x f_h  - \dfrac{\Delta t}{2}\left(\xi\partial_x\partial_t f_h + \xi^2 \partial_x^2 f_h\right)+ \dfrac{\Delta x^2}{6}\xi\partial_x^3 f_h
= \dfrac{1}{\epsilon}Q_h  + O(\Delta x^3, \Delta t^2),
\end{equation}
which shows that the time accuracy of this CLR scheme is of first-order.

As $\Delta t=O(\epsilon^{\alpha})$ and $\Delta x=O(\epsilon^{\beta})$ with $0.5<\alpha,\beta<1$, Equation \eqref{eq:sME-2x} can be reformulated as (with only error terms of leading order),
\begin{equation}
\label{eq:sME-3}
\epsilon \partial_t f_h + \epsilon \xi \partial_x f_h   -\epsilon^{1+\alpha} \dfrac{a}{2}\left(\xi\partial_x\partial_t f_h + \xi^2 \partial_x^2 f_h\right) + \epsilon^{2\beta+1}\dfrac{b}{6}\xi\partial_x^3 f_h  = -\dfrac{1}{\tau}(f_h-f_h^{(eq)}),
\end{equation}
where $a=O(1)$ and $b=O(1)$. From Eq. \eqref{eq:sME-3} we can obtain the first two Chapman-Enskog expansion coefficients of $f_h$
\begin{subequations}
\label{eq:sCE-BGK}
\begin{align}
\epsilon^0:&\qquad f_h^{(0)} = f_h^{(eq)},\\
\epsilon^1:&\qquad  D_0 f_h^{(0)} =   -\dfrac{1}{\tau} f_h^{(1)}
\label{eq:sCE-1}.
\end{align}
\end{subequations}
However, the coefficient $f_h^{(2)}$ cannot be determined due to the appearing of term of $O(\epsilon^{1+\alpha})$. This suggests that CLR scheme is of first-order UP and the Navier-Stokes solution cannot be captured by this scheme under the same numerical resolution as the DUGKS, i.e., $\Delta t=O(\epsilon^{\alpha})$ and $\Delta x=O(\epsilon^{\beta})$ with $0.5<\alpha,\beta<1$.

It is interesting that as $\Delta t=o(\epsilon)$ and $\Delta x=o(\epsilon^{1/2})$, namely $\Delta t=O(\epsilon^{\alpha})$ and $\Delta x=O(\epsilon^{\beta})$ with $\alpha>1$ and $\beta >0.5$, we can obtain the balance equation for $f_h^{(2)}$ exactly. However, this means that the time step $\Delta t$ must resolve the kinetic time scale in order to capture the Navier-Stokes solutions, which indicates that the CLR scheme is not an UP solver and is inapplicable to multiscale flows. The above analysis confirms the fact that it is necessary to consider the collision effect in the reconstruction of the cell-interface distribution function in developing UP kinetic schemes.

\bibliographystyle{elsarticle-num}
\bibliography{references}

\begin{thebibliography}{10}
\expandafter\ifx\csname url\endcsname\relax
  \def\url#1{\texttt{#1}}\fi
\expandafter\ifx\csname urlprefix\endcsname\relax\def\urlprefix{URL }\fi
\expandafter\ifx\csname href\endcsname\relax
  \def\href#1#2{#2} \def\path#1{#1}\fi

\bibitem{ref:Xu-POF2017}
K.~Xu, C.~Liu, A paradigm for modeling and computation of gas dynamics, Phys.
  Fluids 29 (2017) 026101.

\bibitem{ref:GuoBook}
Z.~L. Guo, C.~Shu, Lattice {Boltzmann} method and its applications in
  engineering, World Scientific, Singapore, 2013.

\bibitem{ref:GKS}
K.~Xu, A gas-kinetic {BGK} scheme for the {Navier-Stokes} equations and its
  connection with artificial dissipation and {Godunov} method, J. Comput. Phys.
  171 (2001) 289--335.

\bibitem{ref:semi-Lagrangian}
P.~Santagati, G.~Russo, S.~B. Yun, Convergence of a {semi-Lagrangian} scheme
  for the {BGK} model of the {Boltzmann} equation, SIAM J. Numer. Anal. 50
  (2012) 1111--1135.

\bibitem{ref:IMEX13}
G.~Dimarco, L.~Pareschi, Asymptotic preserving implicit-explicit {Runge-Kutta}
  methods for nonlinear kinetic equations, SIAM J. Numer. Anal. 51 (2013)
  1064--1087.

\bibitem{ref:IMEX17}
G.~Dimarco, L.~Pareschi, Implicit-explicit linear multistep methods for stiff
  kinetic equations, SIAM J. Numer. Anal. 55 (2017) 664--690.

\bibitem{ref:IMEX07}
S.~Pieraccini, G.~Puppo, Implicit-explicit schemes for {BGK} kinetic equations,
  J. Sci. Comput. 32 (2007) 1--28.

\bibitem{ref:IMEX_Hu18}
J.~W. Hu, R.~W. Shu, X.~X. Zhang, Asymptotic-preserving and
  positivity-preserving implicit-explicit schemes for the stiff {BGK} equation,
  SIAM J. Numer. Anal. 56~(2) (2018) 942--973.

\bibitem{ref:UGKS}
K.~Xu, J.~C. Huang, A unified gas-kinetic scheme for continuum and rarefied
  flows, J. Comput. Phys. 229 (2010) 7747--7764.

\bibitem{ref:DUGKS15}
Z.~L. Guo, R.~J. Wang, K.~Xu, Discrete unified gas kinetic scheme for all
  {Knudsen} number flows. {II. Thermal} compressible case, Phys. Rev. E 91
  (2015) 033313.

\bibitem{ref:DUGKS13}
Z.~L. Guo, K.~Xu, R.~J. Wang, Discrete unified gas kinetic scheme for all
  {Knudsen} number flows: {Low-speed} isothermal case, Phys. Rev. E 88 (2013)
  033305.

\bibitem{ref:DUGKS_Rev}
Z.~L. {Guo}, K.~Xu, Progress of discrete unified gas-kinetic scheme for
  multiscale flows, Adv. Aerodyn. 3 (2021) 6.

\bibitem{ref:Review14}
G.~Dimarco, L.~Pareschi, Numerical methods for kinetic equations, Acta Numerica
  23 (2014) 369--520.

\bibitem{ref:AP_Rev_2017}
J.~W. Hu, S.~Jin, Q.~Li, Asymptotic-preserving schemes for multiscale
  hyperbolic and kinetic equations, In R. Abgrall and C.-W. Shu, editors,
  Handbook of Numerical Methods for Hyperbolic Problems, chapter 5, pages
  103-129, North-Holland, 2017.

\bibitem{ref:AP_Rev_2012}
S.~Jin, Asymptotic preserving {(AP)} schemes for multiscale kinetic and
  hyperbolic equations: a review, Riv. Mat. Univ. Parma 3 (2012) 177--216.

\bibitem{ref:Jin1999}
S.~Jin, Efficient asymptotic-preserving ({AP}) schemes for some multiscale
  kinetic equations, SIAM J. Sci. Comput. 21~(2) (1999) 441--454.

\bibitem{ref:Larsen_1983}
E.~W. Larsen, On numerical solutions of transport problems in the diffusion
  limit, Nucl. Sci. Eng. 83~(1) (1983) 90--99.

\bibitem{ref:Larsen_1987}
E.~W. Larsen, J.~Morel, J.~Miller, Asymptotic solutions of numerical transport
  problems in optically thick, diffusive regimes, J. Comput. Phys. 69 (1987)
  283--324.

\bibitem{ref:Klar98}
A.~Klar, An asymptotic-induced scheme for nonstationary transport equations in
  the diffusive limit, SIAM J. Numer. Anal. 35~(3) (1998) 1073--1094.

\bibitem{ref:Jin1995}
S.~Jin, {Runge-Kutta} methods for hyperbolic conservation laws with stiff
  relaxation terms, J. Comput. Phys. 122~(1) (1995) 51--67.

\bibitem{ref:Pareschi2003}
L.~Pareschi, G.~Russo, High order asymptotically strong-stability-preserving
  methods for hyperbolic systems with stiff relaxation, in: Hyperbolic
  Problems: Theory, Numerics, Applications, Springer, 2003, pp. 241--251.

\bibitem{ref:Lowrie02}
R.~B. Lowrie, J.~E. Morel, Methods for hyperbolic systems with stiff
  relaxation, Int. J. Numer. Meth. Fluids 40~(3-4) (2002) 413--423.

\bibitem{ref:Coron91}
F.~Coron, B.~Perthame, Numerical passage from kinetic to fluid equations, SIAM
  J. Numer. Anal. 28~(1) (1991) 26--42.

\bibitem{ref:AP-NS-Klar}
A.~Klar, An asymptotic preserving numerical scheme for kinetic equations in the
  low {Mach} number limit, SIAM J. Numer. Anal. 36 (1999) 1507--1527.

\bibitem{ref:Luc_NS_JCP08}
M.~Bennoune, M.~Lemou, L.~Mieussens, Uniformly stable numerical schemes for the
  {Boltzmann} equation preserving the compressible {Navier--Stokes}
  asymptotics, J. Comput. Phys. 227~(8) (2008) 3781--3803.

\bibitem{ref:AP-NS-FJJCP2010}
F.~Filbet, S.~Jin, A class of asymptotic-preserving schemes for kinetic
  equations and related problems with stiff sources, J. Comput. Phys. 229~(20)
  (2010) 7625--7648.

\bibitem{ref:Luc_RGD2014}
L.~Mieussens, A survey of deterministic solvers for rarefied flows, in: AIP
  Conf. Proc., Vol. 1628, American Institute of Physics, 2014, pp. 943--951.

\bibitem{ref:Tcher01}
F.~G. Tcheremissine, Solution of the {Boltzmann} equation in stiff regime, in:
  Hyperbolic Problems: Theory, Numerics, Applications, Springer, 2001, pp.
  883--890.

\bibitem{ref:KlarJCP99}
A.~Klar, Relaxation scheme for a lattice--{Boltzmann}-type discrete velocity
  model and numerical {Navier--Stokes} limit, J. Comput. Phys. 148~(2) (1999)
  416--432.

\bibitem{ref:AP_NS_Hu2017}
J.~W. Hu, X.~X. Zhang, On a class of implicit-explicit {Runge--Kutta} schemes
  for stiff kinetic equations preserving the {Navier-Stokes} limit, J. Sci.
  Comput. 73~(2-3) (2017) 797--818.

\bibitem{ref:AP_NS_Tong_JCP15}
T.~Xiong, J.~Jang, F.~Y. Li, J.-M. Qiu, High order asymptotic preserving nodal
  discontinuous {Galerkin IMEX }schemes for the {BGK} equation, J. Comput.
  Phys. 284 (2015) 70--94.

\bibitem{ref:Klar_CF2006}
M.~Sea{\"\i}d, A.~Klar, Asymptotic-preserving schemes for unsteady flow
  simulations, Comput. Fluids 35~(8-9) (2006) 872--878.

\bibitem{AP-NS-Bos2017}
S.~Boscarino, L.~Pareschi, On the asymptotic properties of {IMEX Runge--Kutta}
  schemes for hyperbolic balance laws, J. Comput. Appl. Math. 316 (2017)
  60--73.

\bibitem{ref:Jin_Levermore1996}
S.~Jin, C.~D. Levermore, Numerical schemes for hyperbolic conservation laws
  with stiff relaxation terms, J. Comput. Phys. 126~(2) (1996) 449--467.

\bibitem{ref:Li-TS-2019}
J.~Q. Li, Two-stage fourth order: temporal-spatial coupling in computational
  fluid dynamics {(CFD)}, Adv. Aerodyn. 1 (2019) 3.

\bibitem{ref:Li-Du-2016}
J.~Q. Li, Z.~F. Du, A two-stage fourth order time-accurate discretization for
  {Lax-Wendroff} type flow solvers {I. Hyperbolic} conservation laws, SIAM J.
  Sci. Comput. 38 (2016) A3046--A3069.

\bibitem{ref:Pan-Xu-Li-Li-2016}
L.~Pan, K.~Xu, Q.~B. Li, J.~Q. Li, An efficient and accurate two-stage
  fourth-order gas-kinetic scheme for the {Euler} and {Navier-Stokes}
  equations, J. Comput. Phys. 326 (2016) 197--221.

\bibitem{ref:ChenXu-2015}
S.~Z. Chen, K.~Xu, A comparative study of an asymptotic preserving scheme and
  unified gas-kinetic scheme in continuum flow limit, J. Comput. Phys. 288
  (2015) 52--65.

\bibitem{ref:GDVM-WP}
P.~Wang, M.~T. Ho, L.~Wu, Z.~L. Guo, Y.~H. Zhang, A comparative study of
  discrete velocity methods for low-speed rarefied gas flows, Comput. Fluids
  161 (2018) 33--46.

\bibitem{ref:ME}
R.~F. Warming, B.~J. Hyett, The modified equation approach to the stability and
  accuracy analysis of finite-difference methods, J. Comput. Phys. 14 (1974)
  159--179.

\bibitem{ref:Li-2013}
J.~Q. Li, Z.~C. Yang, The von {Neumann} analysis and modified equation approach
  for finite difference schemes, Appl. Math. Comput. 225 (2013) 610--621.

\bibitem{ref:ChapmanBook}
S.~Chapman, T.~Cowling, The mathematical theory of non-uniform gases, Cambridge
  University Press, 1970.

\bibitem{ref:CercignaniBook}
C.~Cercignani, The {Boltzmann} equation and its applications, Springer-Verlag,
  New York, 1988.

\bibitem{ref:Gorban2014}
A.~N. Gorban, I.~V. Karlin, Hilbert's 6th problem: exact and approximate
  hydrodynamic manifolds for kinetic equations, Bull. Amer. Math. Soc. 51
  (2014) 187--246.

\bibitem{ref:RegCE_Rose1989}
P.~Rosenau, Extending hydrodynamics via the regularization of the
  {Chapman-Enskog} expansion, Phys. Rev. A 40~(12) (1989) 7193--7196.

\bibitem{ref:NormCE_Slemrod1989}
M.~Slemrod, A renormalization method for the {Chapman-Enskog} expansion that
  yields truncation stability, Phys. D 109 (1997) 257--273.

\bibitem{ref:SumCE_Gorban1982}
A.~N. Gorban, I.~V. Karlin, Structure and approximations of the
  {Chapman-Enskog} expansion for the linearized {Grad} equations, Transp.
  Theor. Stat. Phys. 21 (1992) 101--117.

\bibitem{ref:NormCE_Slemrod2012}
M.~Slemrod, {Chapman-Enskog} ${\Rightarrow}$ viscosity-capillarity, Quart.
  Appl. Math. 70~(3) (2012) 613--624.

\bibitem{ref:LB_Wagner2006}
A.~J. Wagner, Thermodynamic consistency of liquid-gas lattice {Boltzmann}
  simulations, Phys. Rev. E 74~(5) (2006) 056703.

\bibitem{ref:LB_HuangRZ2016}
R.~Z. Huang, H.~Y. Wu, Third-order analysis of pseudopotential lattice
  {Boltzmann} model for multiphase flow, J. Comput. Phys. 327 (2016) 121--139.

\bibitem{ref:LB_ZhengL2017}
L.~{Zheng}, Q.~L. {Zhai}, S.~{Zheng}, Analysis of force treatment in the
  pseudopotential lattice {Boltzmann} equation method., Phys. Rev. E 95~(4)
  (2017) 043301.

\bibitem{ref:BGK_Burnett}
K.~{Xu}, Z.~H. {Li}, Microchannel flow in the slip regime: gas-kinetic
  {BGK-Burnett} solutions, J. Fluid Mech. 513 (2004) 87--110.

\bibitem{ref:Junk-Yong2002}
M.~Junk, W.-A. Yong, Rigorous {Navier-Stokes} limit of the lattice {Boltzmann}
  equation, Asymptotic Anal. 35 (2002) 165--185.

\bibitem{ref:Junk2005}
M.~{Junk}, A.~{Klar}, L.-S. {Luo}, Asymptotic analysis of the lattice
  {Boltzmann} equation, J. Comput. Phys. 210~(2) (2005) 676--704.

\bibitem{ref:DUGKS-A3}
L.~H. Zhu, S.~Z. Chen, Z.~L. Guo, {dugksFoam}: An open source {OpenFOAM} solver
  for the {Boltzmann} model equation, Comput. Phys. Commun. 213 (2017)
  155--164.

\bibitem{ref:DUGKS-A2}
L.~H. Zhu, X.~F. Yang, Z.~L. Guo, Numerical study of nonequilibrium gas flow in
  a microchannel with a ratchet surface, Phys. Rev. E 95 (2017) 023113.

\bibitem{ref:DUGKS-A1}
L.~H. Zhu, X.~F. Yang, Z.~L. Guo, Thermally induced rarefied gas flow in a
  three-dimensional enclosure with square cross-section, Phys. Rev. Fluids 2
  (2017) 123402.

\bibitem{ref:ESWENO}
N.~K. Yamaleev, M.~H. Carpenter, Third-order energy stable {WENO} scheme, J.
  Comput. Phys. 228~(8) (2013) 3025--3047.

\bibitem{ref:BirdBook}
G.~A. Bird, J.~M. Brady, Molecular gas dynamics and the direct simulation of
  gas flows, Vol.~42, Clarendon press Oxford, 1994.

\bibitem{ref:JiangDW}
D.~W. Jiang, M.~L. Mao, J.~Li, X.~G. Deng, An implicit parallel {UGKS} solver
  for flows covering various regimes, Adv. Aerodyn. 1 (2019) 8.

\bibitem{ref:HarrisBook}
S.~Harris, An introduction to the theory of the {Boltzmann} equation, Dover
  Publications, 2011.

\bibitem{ref:Struchtrup2004}
H.~Struchtrup, Stable transport equations for rarefied gases at high orders in
  the {Knudsen} number, Phys. Fluids 16 (2004) 3921--3934.

\bibitem{ref:Kauf2010scale}
P.~{Kauf}, M.~{Torrilhon}, M.~{Junk}, Scale-induced closure for approximations
  of kinetic equations, J. Stat. Phys. 141~(5) (2010) 848--888.

\end{thebibliography}

\end{document}